\documentclass{amsart}
\usepackage{amsmath,amsthm,amsfonts,amssymb,latexsym}

\newtheorem{theorem}{Theorem}[section]
\newtheorem{proposition}[theorem]{Proposition}
\newtheorem{lemma}[theorem]{Lemma}

\newtheorem{corollary}[theorem]{Corollary}
\newtheorem{observation}[theorem]{Observation}
\theoremstyle{definition}

\newtheorem{remark}[theorem]{Remark}

\newcommand{\U}{\mathcal U}
\newcommand{\V}{\mathcal V}
\newcommand{\W}{\mathcal W}
\newcommand{\F}{\mathcal F}
\newcommand{\N}{\mathcal N}
\newcommand{\w}{\omega}
\newcommand{\SF}{\mathsf{SF}}
\newcommand{\UF}{\mathsf{UF}}
\newcommand{\da}{\downarrow}

\newcommand{\cf}{\mathrm{cf}}
\newcommand{\e}{\varepsilon}
\newcommand{\IR}{\mathbb R}

\newcommand{\Ra}{\Rightarrow}
\newcommand{\uss}{\underset{\leftarrow}{\subset}}

\newcommand{\la}{\langle}
\newcommand{\obr}{\mathrm{rng}\,}
\newcommand{\ra}{\rangle}
\newcommand{\uhr}{{\upharpoonright}}

\newcommand{\id}{\mathrm{id}}
\newcommand{\LL}{\mathcal L}

\newcommand{\cen}{\mathsf{cen}}

\title[$o$-Boundedness of free topological groups]{$o$-Boundedness of free topological groups}

\author[T. Banakh,  D. Repov\v{s},  and L. Zdomskyy]{Taras Banakh, Du\v{s}an Repov\v{s},  and Lyubomyr Zdomskyy}

\address{Department of Mathematics, Ivan Franko National University of Lviv, Ukraine and\\
Instytut Matematyki, Uniwersytet Humanistyczno-Przyrodniczy Jana Kochanowskiego w Kielcach, Poland.}
\email{tbanakh@yahoo.com}

\address{Faculty of Mathematics and Physics, and Faculty of Education,
     University of Ljubljana, Jadranska 19,
Ljubljana, Slovenija 1000.}
\email{dusan.repovs@guest.arnes.si}

\address{Kurt G\"odel Research Center for Mathematical Logic,
University of Vienna, W\"ahringer Stra\ss e 25, A-1090 Wien, Austria.}
\email{lzdomsky@gmail.com}

\thanks{The  authors were supported in part by the Slovenian
Research Agency grants P1-0292-0101-04, J1-9643-0101, and
BI-UA/07-08-001. The third author would also like to acknowledge the support
of the FWF grant
P19898-N18.}

\begin{document}
\begin{abstract}
Assuming the absence of Q-points (which is consistent with ZFC) we prove that
the free topological group $F(X)$ over a Tychonov
space $X$ is $o$-bounded if and only if every continuous metrizable
image $T$ of $X$ satisfies the selection principle
$\bigcup_{\mathit{fin}}(\mathcal O,\Omega)$ (the latter means that
for every sequence $\la u_n\ra_{n\in\w}$ of open covers of $T$ there
exists a sequence $\la v_n\ra_{n\in\w}$ such that $v_n\in
[u_n]^{<\w}$ and for every  $F\in [X]^{<\w}$ there exists $n\in\w$
with $F\subset\cup v_n$).
 This characterization gives a consistent answer to a problem posed by C.~Hernandes,
D.~Robbie,  and M.~Tkachenko in 2000.
\end{abstract}
\maketitle

\section{Introduction}

In this paper we present a (consistent) characterization of topological spaces $X$ whose free topological group $F(X)$ is $o$-bounded, thus resolving the corresponding problem posed by
 C.~Hernandes, D.~Robbie,  and M.~Tkachenko in \cite{HRT00}.

We recall that a topological group $G$ is said to be {\em $o$-bounded\/}
if for every sequence $\la U_n\ra_{n\in\w}$ of neighborhoods of the neutral element $e$ of $G$ there is a sequence
$\la F_n\ra_{n\in\w}$ of finite subsets of $G$ such that
$G=\bigcup_{n\in\w}U_n\cdot F_n$.
 This notion was
introduced by O.~Okunev and M.~Tkachenko as a covering counterpart of a $\sigma$-bounded group,
see e.g. \cite{Tk98} for the discussion of this subject.
It is clear that each $\sigma$-bounded group (that is, a subgroup of a $\sigma$-compact topological group) is
$o$-bounded.

Our aim is to detect topological spaces $X$ with $o$-bounded free (abelian) topological group. By a free abelian topological group over a Tychonov space $X$ we understand an abelian topological group $A(X)$ that contains $X$ as a subspace and such that each continuous map  $\phi:X\to G$ into an abelian topological group $G$ extends to a  unique continuous group homomorphism $\Phi:A(X)\to G$. Deleting the adjective ``abelian'', we obtain the definition of a free topological group $F(X)$ over $X$.

Tychonov spaces $X$ with $o$-bounded free abelian topological group $A(X)$ were characterized by the third author \cite{Zd06} as spaces  all of whose continuous metrizable images are Scheepers spaces. We recall that a topological space $X$ satisfies the selection principle $\bigcup_{\mathit{fin}}(\mathcal O,\mathcal O)$ (resp. $\bigcup_{\mathit{fin}}(\mathcal O,\Omega)$) or else is said to be {\em Menger} (resp. {\em Scheepers}) if for any sequence $\la\mathcal U_n\ra_{n\in\w}$ of open covers of $X$
each cover $\U_n$ contains a finite subfamily $\V_n\subset\U_n$ such that
$\la \cup\V_n\ra_{n\in\w}$ is a cover (an $\w$-cover) of $X$. A cover $\U$ of a set $X$ is called an {\em $\w$-cover} if each finite subset $F\subset X$ lies in some element $U\in \U$ of the cover $\U$.

It is clear that each Scheepers space is a Menger space. The converse follows from the inequality $\mathfrak u<\mathfrak g$ (which holds in some models of ZFC), see \cite{Zd05}. On the other hand, under CH there is a Menger metrizable space that fails to be Scheepers, see \cite[Theorem~2.8]{JMSS}.

The following theorem proved in \cite{Zd06} characterizes Tychonov spaces $X$ with $o$-bounded free abelian topological group $A(X)$.

\begin{theorem}\label{ab_ch} The free abelian topological group $A(X)$ over a Tychonov space $X$ is $o$-bounded if and only if each continuous metrizable image of $X$ is Scheepers.
\end{theorem}

In this paper, assuming the absence of $Q$-points, we shall prove a
similar characterization of $o$-bounded free topological groups
$F(X)$. By a {\em $Q$-point\/} we understand a free ultrafilter $\U$
on $\w$ such that each finite-to-one function $f:\w\to\w$ is
injective on a suitable set $U\in\U$.

The following characterization  is the principal result of this paper\footnote{It was announced in \cite{BZ06} that Theorem~\ref{main} is true without
any additional set-theoretic assumptions, but that proof does
 not work under the existence of a $Q$-point, see our Theorem~\ref{ess_q}.}.

\begin{theorem}\label{main} If there is no $Q$-point, then for every Tychonov space $X$ the following conditions are equivalent:
\begin{enumerate}
\item All finite powers of the free topological group $F(X)$ are
$o$-bounded;
\item The free topological group $F(X)$ is $o$-bounded;
\item The free abelian topological group $A(X)$ is $o$-bounded;
\item Every continuous metrizable image of $X$ is Scheepers.
\end{enumerate}
\end{theorem}

The existence of $Q$-points is independent from ZFC. One can easily
construct $Q$-points under Continuum Hypothesis. On the other hand,
since under NCF for every ultrafilter $\F$ there exists a monotone
surjection $\phi:\w\to\w$ such that $\phi(\F)$ is generated by
$\mathfrak u$ sets and $\mathfrak u<\mathfrak d$ by
\cite[Theorem~1]{afterall}, and obviously no $Q$-point can be
generated by less than $\mathfrak d$ sets, NCF implies that there
are no $Q$-points. In particular, there are no $Q$-points in the
models of $\mathfrak u<\mathfrak g$\footnote{The assertion
``$\mathfrak b=\mathfrak d$ and there are no $Q$-points'' is
consistent as well, see \cite{Mi80}.}.

 At this point we would like to mention that one
of the most important problems regarding the $o$-boundedness is
 whether this property is preserved by squares of topological groups. This
 problem was first solved in \cite{MST} under additional
 set-theoretic assumptions, which were further weakened to $\mathfrak r\geq\mathfrak d$
 by Mildenberger \cite{Mi08}. The apparently most intriguing part of
 this problem left open after \cite{MST} and \cite{Mi08} is whether
 the square of an $o$-bounded group is $o$-bounded provided NCF
 (or even $\mathfrak u<\mathfrak g$) holds. From the above it
 follows  that  free topological groups can not be used to solve
 this problem.

The proof of Theorem~\ref{main} is nontrivial and relies on the theory of multicovered spaces developed in \cite{BZ06}. The necessary information related to multicovered spaces is collected in Section~\ref{mult_sp}. In that section we also define the operation of the semi-direct product of multicovered spaces. In Section~\ref{f_meng} we recall some information related to the $\mathsf F$-Menger property in multicovered spaces and prove an important Theorem~\ref{semiprod} on preservation of the $[\F]$-Menger property by semi-direct products of multicovered spaces. In Section~\ref{sem_mult_gr} we characterize the
$\mathsf F$-Menger property in (semi)multicovered groupoids.
Those are multicovered spaces endowed with a binary and unary operations that are compactible with the multicover in a suitable sense. In
Sections~\ref{F-Menger-groups}, \ref{top-monoids}, \ref{top_mon_s} we apply the obtained results
about semi-multicovered groupoids to detecting the $\mathsf F$-Menger property in topological groups and topological monoids. These results provide us with tools for the proof of Theorem~\ref{main} which is given in Section~\ref{letztlich}.

\section{Multicovered spaces} \label{mult_sp}

By a {\em multicovered space} we understand a pair $(X,\mu)$ consisting of a set $X$ and a family $\mu$ of covers of $X$. Such a family $\mu$ is called the {\em multicover} of $X$. Multicovered spaces naturally appear in many situations. In particular,
\begin{itemize}
\item each topological space $X$ has the canonical multicover $\mu_{\mathcal O}$ consisting of all open covers of $X$;
\item each uniform space $(X,\U)$ possesses the canonical multicover\newline
$\mu_\U=\big\{\{U(x):x\in X\}:U\in\U\big\}$ consisting of uniform covers of $X$;
\item each metric space $(X,d)$ has the canonical multicover\newline $\mu_d=\big\{\{B(x,\e):x\in X\}:\e>0\big\}$ consisting of covers by balls of fixed radius;
\item each topological group $G$ with topology $\tau$ carries four natural multicovers:
\begin{itemize}
\item the left multicover $\mu_L=\big\{\{xU:x\in G\}:1\in U\in\tau\big\}$;
\item the right multicover  $\mu_R=\big\{\{Ux:x\in G\}:1\in U\in\tau\big\}$;
\item the two-sided multicover  $\mu_{L\wedge R}=\big\{\{Ux\cap xU:x\in G\}:1\in U\in\tau\big\}$;
\item the R\"olke multicover  $\mu_{L\vee R}=\big\{\{UxU:x\in G\}:1\in U\in\tau\big\}$.
\end{itemize}
\end{itemize}

A multicovered space $(X,\mu)$ is  {\em centered} if for any covers $u,v\in\mu$ there
is a cover $w\in\mu$ such that each $w$-bounded subset $B\subset X$ is both $u$-bounded and $v$-bounded. We define a subset $B\subset X$ to be {\em $u$-bounded} with respect to a cover $u$ of $X$  if $B\subset \cup u'$ for some finite subfamily $u'\subset u$.

Observe that all the examples of multicovered spaces presented above are centered.
Each multicover $\mu$ on $X$ can be transformed into a centered multicover $\cen(\mu)$ consisting of the covers $$u_1\wedge\dots\wedge u_n=\{U_1\cap\dots \cap U_n:\forall i\le n\;\; U_i\in u_i\}\mbox{ where $u_1,\dots,u_n\in\mu$}.$$
From now on {\em all multicovered spaces will be assumed to be centered.}

A multicovered space $X$ is called {\em $\w$-bounded} if each cover $u\in\mu$ contains a countable subcover $u'\subset u$. If, in addition, there is a cover $v\in\mu$ such that each $v$-bounded subset of $X$ is $u'$-bounded, then the multicovered space $(X,\mu)$ is {\em properly $\w$-bounded}.

By the {\em cofinality} $\cf(X,\mu)$ of a multicovered space $(X,\mu)$ we understand  the smallest cardinality $|\mu'|$ of a subfamily $\mu'\subset\mu$ such that for every cover $u\in\mu$ there is a cover $u'\in\mu'$ such that each $u'$-bounded subset $B\subset X$ is $u$-bounded. For example, the multicover of any metric space has countable cofinality.

Multicovered spaces are objects of a category whose morphisms are uniformly bounded maps. We define a map $f:X\to Y$ between two multicovered spaces $(X,\mu_X)$ and $(Y,\mu_Y)$ to be {\em uniformly
 bounded} if for every cover $v\in\mu_Y$ there
is a cover $u\in\mu_X$ such that for every $u$-bounded subset
 $B\subset X$ the image $f(B)$ is $v$-bounded in $Y$.
Observe that every
uniformly continuous  map $f:X\to Y$ between uniform spaces is uniformly bounded
with respect to the multicovers induced by the
 uniformities on $X$ and $Y$, respectively.

Two multicovers $\mu_0,\mu_1$ on a set $X$ are {\em equivalent\/} if the identity maps $\mathrm{id}:(X,\mu_0)\to (X,\mu_1)$ and $\mathrm{id}:(X,\mu_1)\to(X,\mu_0)$ are uniformly bounded (and hence are isomorphisms in the category of multicovered spaces).

\begin{proposition} \label{pr2_1} A multicovered space $(X,\mu)$ is properly $\w$-bounded if and only if $\mu$ is equivalent to a multicover $\nu$ on $X$ consisting of countable disjoint covers.
\end{proposition}

\begin{proof}
Suppose that $\mu$ is an equivalent to
a multicover $\nu$ of $X$ consisting of countable covers. Given any $u\in\mu$, we can
find $v\in\nu$ such that each $v$-bounded subset of $X$ is $u$-bounded.
This means that for every $V\in v$ there exists $u_V\in [u]^{<\w}$ such  that
$V\subset\cup u_V$. Set $u'=\bigcup_{V\in v}u_V$ and observe that $u'$
is a countable subcover of $u$ such that each $v$-bounded subset is $u'$-bounded.

Let $u_1\in\mu$ be such that each $u_1$-bounded subset of $X$ is $v$-bounded.
Then each $u_1$-bounded subset of $X$ is $u'$-bounded, and hence $\mu$ is
properly $\w$-bounded.
\smallskip

Now, suppose that $\mu$ is properly $\w$-bounded. For every $u\in\mu$
fix $u'\in [u]^\w$ and $u''\in\mu$ such that $\cup u'=X$ and each $u''$-bounded
subset of $X$ is $u'$-bounded. A direct verification shows that the multicover
 $\mu'=\{u':u\in\mu\}$ of $X$ is equivalent to $\mu$.
Let us write $u'$ in the form $\{U_n:n\in\w\}$ and let
 $V_n=U_{n}\setminus\bigcup_{i< n}V_n$, $n\in\w$, and $v_u=\{V_n:n\in\w\}$.
 Then the multicover $\nu=\{v_u:u\in\mu\}$ of $X$ is equivalent to $\mu'$
(and hence to $\mu$) and consists of countable disjoint covers.
\end{proof}

Sometimes we shall need the following obvious observation:

\begin{observation} \label{trivial}
If (centered) multicovers $\mu_0$ and $\mu_1$ of a set $X$ are equivalent, then for every $\nu_0\subset\mu_0$
there exist (centered) submulticovers $\nu'_0\subset \mu_0$ and $\nu_1\subset\mu_1$ equivalent to $\nu'_0$ such that
$\nu_0\subset\nu'_0$ and
$|\nu_1|\leq |\nu'_0|\le\max\{|\nu_0|,\w\}$.
\end{observation}

A multicovered space $(X,\mu)$ is {\em uniformizable} if $\mu$ is equivalent to the multicover $$\mu_\U=\big\{\{U(x):x\in X\}:U\in\U\big\}$$ generated by a suitable uniformity $\U$ on $X$ (here for an entourage $U\in\U$ by $U(x)=\{y\in X:(x,y)\in U\}$ we denote the $U$-ball centered at $x$). There is a simple criterion of the uniformizability of $\w$-bounded multicovered spaces.

\begin{proposition} \label{unif} An $\w$-bounded multicovered space is uniformizable
if and only if it is centered and properly $\w$-bounded.
\end{proposition}

\begin{proof} To prove the ``only if'' part, assume that the multicover $\mu$ is equivalent to the multicover $\mu_\U$ generated  by some uniformity $\U$ on $X$.
Since the multicover $\mu_\U$ is centered, so is the multicover $\mu$. Next, we check that $\mu$ is properly $\w$-bounded.

 Given any $u\in\mu$, we can find an entourage $U\in\U$ such that
each $U$-ball $U(x)=\{x'\in X:(x,x')\in U\}$ is $u$-bounded. Let $V\in\U$ be an entourage such that $V=V^{-1}$ and $V\circ V\subset U$.

The $\w$-boundedness of the multicover $\mu$ implies that of $\mu_\U$. Consequently, there is a  countable subset $C\subset X$ such that $\bigcup_{x\in C}V(x)=X$. By the choice of the entourage $U$, there is a countable subfamily $u'\subset u$ such that each $U$-ball $U(c)$, $c\in C$, is $u'$-bounded. Observe that for every $x\in X$ we can find $c\in C$ with $x\in V(c)$. Consequently, $V(x)\subset V\circ V(c)\subset U(c)\subset \cup u'_f$ for some finite subfamily $u'_f\subset u'$. This means that each $V$-ball $V(x)$, $x\in X$, is $u'$-bounded.

Since the multicovers $\mu$ and $\mu_\U$ are equivalent, we can find a cover $v\in\mu$  such that each $v$-bounded subset of $X$ is $\{V(x)\}_{x\in X}$-bounded and hence $u'$-bounded, witnessing that the multicovered space $(X,\mu)$ is properly $\w$-bounded.
\medskip

To prove the ``if'' part, assume that a multicovered space $(X,\mu)$ is centered and properly $\w$-bounded. By Proposition~\ref{pr2_1}, we can assume that the multicover $\mu$ consists of countable disjoint covers. For every finite subfamily $\mu'\subset\mu$ consider the entourage of the diagonal
$$W_{\mu'}=\{(x,y)\in X\times X:\forall u\in\mu'\;\exists U\in u\;\;x,y\in U\}.$$
Those entourages generate a uniformity $\W$ on $X$ such that the multicover $\mu_\W$ consisting of $\W$-uniform covers of $X$ is equivalent to the multicover $\mu$.
\end{proof}

Now we describe some operations in the category of multicovered spaces.

Each subset $A$ of a multicovered space $(X,\mu)$ carries the {\em induced multicover} $\mu\uhr A$ consisting of the covers $u\uhr A=\{U\cap A:U\in u\}$, $u\in\mu$. The multicovered space $(A,\mu\uhr A)$ is called a {\em subspace} of $(X,\mu)$.

Each map $f:X\to Y$ from a set $X$ to a multicovered space $(Y,\mu_Y)$ induces the multicover $f^{-1}(\mu)$ consisting of the covers $f^{-1}(u)=\{f^{-1}(U):U\in u\}$, $u\in\mu_Y$.

For two multicovers $\mu_0,\mu_1$ on a set $X$ their {\em meet} is the multicover
$$\mu_0\wedge\mu_1=\{u\wedge v:u\in\mu_0,\;v\in\mu_1\}$$ where
$$u\wedge v=\{U\cap V:U\in u,\; V\in v\}.$$

The product $X\times Y$ of two multicovered spaces $(X,\mu_X)$ and $(Y,\mu_Y)$ possesses the multicover
$$\mu_X\boxtimes\mu_Y=\{u\boxtimes v:u\in\mu_X,\;v\in\mu_Y\}\mbox{ where }u\boxtimes v=\{U\times V:U\in u,\; V\in v\}.$$
The obtained multicovered space $(X\times Y,\mu_X\boxtimes\mu_Y)$ is called the {\em direct product} of the multicovered spaces $(X,\mu_X)$ and $(Y,\mu_Y)$.

Besides the multicover $\mu_X\boxtimes\mu_Y$ on the Cartesian product $X\times Y$ of multicovered spaces $(X,\mu_X)$ and $(Y,\mu_Y)$ there are three less evident multicovers:
\begin{itemize}
\item  $\mu_X{\ltimes}\mu_Y=\big\{\{U\times V:U\in u,\; V\in v_U\}:u\in\mu_X,\;\{v_U\}_{U\in u}\subset\mu_Y\big\}$, the multicover of the left semi-direct product,
\item $\mu_X{\rtimes}\mu_Y=\big\{\{U\times V:V\in v,\; U\in u_V\}:v\in\mu_Y,\;\{u_V\}_{V\in v}\subset\mu_X\big\}$, the multicover of the right semi-direct product, and
\item $\mu_X{\bowtie}\mu_Y=(\mu_X{\ltimes}\mu_Y)\wedge(\mu_X{\rtimes}\mu_Y)$, the multicover of the semi-direct product.
\end{itemize}

The multicovered spaces
$$X\ltimes Y=(X\times Y,\mu_X{\ltimes}\mu_Y),\;\; X\rtimes Y=(X\times Y,\mu_X{\rtimes}\mu_Y),\;\;
X\bowtie Y=(X\times Y,\mu_X{\bowtie}\mu_Y),$$
are called respectively: the {\em left semi-direct product}, the {\em right semi-direct product},
and the {\em semi-direct product\/} of the multicovered spaces $(X,\mu_X)$ and $(Y,\mu_Y)$.

The following proposition shows that the (semi)direct products respect the equivalence relation.

\begin{proposition} \label{ochev}
Let $\la \mu_X,\nu_X \ra$ and $\la \mu_Y,\nu_Y \ra$ be pairs of equivalent centered multicovers
of sets $X$ and $Y$, respectively. Then $\mu_X\boxtimes\mu_Y$, $\mu_X\ltimes\mu_Y$,
$\mu_X\rtimes\mu_Y$, and $\mu_X\bowtie\mu_Y$ are equivalent to $\nu_X\boxtimes\nu_Y$,
 $\nu_X\ltimes\nu_Y$, $\nu_X\rtimes\nu_Y$, and $\nu_X\bowtie\nu_Y$, respectively.
\end{proposition}

\begin{proof} We shall prove the equivalence of the multicovers
$\mu_X\ltimes\mu_Y$ and $\nu_X\ltimes\nu_Y$. For the other pairs of multicovers the proof is analogous.

In order to prove the uniform boundedness of the identity map $$(X\times Y,\mu_X\ltimes\mu_Y)\to(X\times Y,\nu_X\ltimes\nu_Y),$$ take any cover $w\in \nu_X\ltimes\nu_Y$ and find a cover $u\in\nu_X$ and a family of covers $\{v_U\}_{U\in u}\subset \nu_Y$ such that $w=\{U\times V:U\in u,\;V\in v_U\}$.

Since the cover $\mu_X$ is equivalent to $\nu_X$, there is a cover $u'\in\mu_X$ such that each $u'$-bounded subset of $X$ is $u$-bounded. Consequently, for each $U'\in u'$ there is a finite subfamily $u_{U'}\subset u$ with $U'\subset \cup u_{U'}$. Since the multicover $\mu_Y$ is equivalent to the centered multicover $\nu_Y$, there is a
cover $v'_{U'}\in\mu_Y$ such that each $v'_{U'}$-bounded subset of $Y$ is $v_U$-bounded for every set $U\in u'_{U'}$.
Now consider the cover $$w'=\{U'\times V':U'\in u',\;\V'\in v'_{U'}\}\in\mu_X\ltimes \mu_Y$$ and observe that each $w'$-bounded subset of $X\times Y$ is $w$-bounded.

Analogously we prove the uniform boundedness of the identity map $$(X\times Y,\nu_X\ltimes\nu_Y)\to(X\times Y,\mu_X\ltimes\mu_Y).$$
\end{proof}

The (semi-)direct product operations also respect the (proper) $\w$-boundedness.
This follows from Propositions~\ref{pr2_1}, \ref{ochev} and Observation~\ref{trivial} above.

\begin{proposition}\label{pob} Let $(X,\mu_X)$ and $(Y,\mu_Y)$ be two multicovered spaces. If the multicovers $\mu_X$ and $\mu_Y$ are (properly) $\w$-bounded, then so are the multicovers \mbox{$\mu_X\boxtimes\mu_Y$,} $\mu_X\ltimes\mu_Y$, $\mu_X\rtimes\mu_Y$, and $\mu_X\bowtie\mu_Y$ on $X\times Y$.
\end{proposition}

Finally we prove the following useful reduction lemma.

\begin{lemma} \label{to_count} Let $(X,\mu_X)$ and $(Y,\mu_Y)$ be two properly $\w$-bounded multicovered \mbox{spaces}. For every countable subfamily $\mu'\subset\mu_X{\bowtie}\mu_Y$ there are countable subfamilies $\mu'_X\subset\mu_X$ and $\mu'_Y\subset \mu_Y$ such that the identity map $\id:(X\times Y,\mu'_X\bowtie\mu'_Y)\to(X\times Y,\mu')$ is uniformly bounded.
\end{lemma}

\begin{proof}
Observation~\ref{trivial} and Propositions~\ref{pr2_1}, \ref{ochev} reduce
the proof of the lemma to the case when $\mu_X$
and $\mu_Y$ consist of countable disjoint covers. Also, there is no loss of generality in assuming that
$\mu'=\{w\}$, where $w=w_0\wedge w_1$ for some $w_0\in\mu_X\rtimes\mu_Y$ and
$w_1\in\mu_X\ltimes\mu_Y$. In its turn, we can write $w_0$ and $w_1$ in the form
$w_0=\{U\times V:V\in v,\; U\in u_V\}$, where $v\in\mu_Y$ and $\la u_V\ra_{V\in v}\in\mu_X^v$;
and $w_1=\{U\times V:U\in u,\; V\in v_U\}$, where $ u\in\mu_X$ and $\la v_U\ra_{U\in u}\in\mu_Y^u$.
Now, any countable centered $\mu'_X\subset\mu_X$ and $\mu'_Y\subset\mu_Y$
containing  $\{u\}\cup\{u_V:V\in v\}$ and $\{v\}\cup\{v_U:U\in u\}$, respectively,
 are easily seen to be as asserted.
\end{proof}


\section{$\mathsf F$-Menger multicovered spaces} \label{f_meng}

The notions of a Menger or Scheepers topological space can be easily generalized to multicovered spaces.
Namely, we define a multicovered space $(X,\mu)$ to be
{\em Menger} (resp. {\em Scheepers}) if for any sequence of covers $\la u_n\ra_{n\in\w}\in\mu^\w$ there is a cover (resp. $\w$-cover) $\la B_n\ra_{n\in\w}$ of $X$ by $u_n$-bounded sets $B_n\subset X$, $n\in\w$.
These two properties are particular cases of the $\mathsf F$-Menger property
of a multicovered space, where $\mathsf F$ is a suitable family of semifilters.

By a {\em semifilter} we understand a family $\F$ of infinite subsets of $\w$ such that for each set
 $A\in\F$ the free filter $\F_A=\{B\subset\w:|A\setminus B|<\aleph_0\}$ generated by
$A$ lies in
 $\F$. Intuitively, sets belonging to a given semifilter
can be thought of  as large in a suitable sense.
 An evident example of a semifilter is any free
(ultra)filer on $\w$. By $\SF$ we denote the family
 of all semifilters and by $\UF$ the subfamily of $\SF$ consisting
of all free ultrafilters (so,
 $\UF=\beta\w\setminus\w$).
More information on semifilters can be found in \cite{BZd} and \cite{BZs1}.

There is an important equivalence relation on the set of semifilters $\SF$, called the {\em coherence relation}, which is defined with help of finite-to-finite multifunctions.

By a multifunction on $\w$ we understand a subset $\Phi\subset\w\times\w$ thought of as a set-valued function $\Phi:\w\Ra\w$ assigning to each number $n\in\w$ the subset $\Phi(n)=\{m\in\w:(n,m)\in\Phi\}$ of $\w$. Each multifunction $\Phi$ has the inverse $\Phi^{-1}=\{(n,m):(m,n)\in\Phi\}$. A multifunction $\Phi:\w\Ra\w$ is called {\em finite-to-finite} if for every  $a\in \w$ the sets $\Phi(a)$ and $\Phi^{-1}(a)$ are finite and non-empty.

We say that a semifilter $\F$ is {\em subcoherent} to a semifilter $\U$ and denote this by $\F\Subset\U$ if there is a finite-to-finite multifunction $\Phi:\w\Ra\w$ such that $\Phi(\F)=\{\Phi(F):F\in\F\}\subset\U$. A semifilters $\F,\U$ on $\w$ are {\em coherent} (denoted by $\F\asymp \U$) if $\F\Subset\U$ and $\U\Subset\F$.
The coherence of semifilters is an equivalence relation dividing the set $\SF$ of
semifilters into the coherence classes $[\F]=\{\U\in\SF:\U\asymp\F\}$ of semifilters $\F$. By Proposition 5.5.2 and Theorem 5.5.3 of \cite{BZd}, a semifilter $\F$ is coherent to a filter $\U$ if and only if $\varphi(\F)=\varphi(\U)$ for some monotone surjection $\varphi:\w\to\w$. It is consistent that all free ultrafilters on $\w$ are coherent, see \cite{BS87}.
 This statement, referred to as NCF (the Near Coherence of Filters)
follows from $\mathfrak u<\mathfrak g$ but contradicts the Martin Axiom,
see \cite{BL89}.

Let $\mathsf F\subset\SF$ be a family of semifilters. Following \cite{BZ06}, we define an indexed cover $\la U_n\ra_{n\in\w}$ of a set $X$ to be an {\em $\mathsf F$-cover} if there is a semifilter $\F\in\mathsf F$ such that for every $x\in X$ the set $\{n\in\w:x\in U_n\}$ belongs to $\F$. Observe that $\la U_n\ra_{n\in\w}$ is an $\w$-cover if and only if it is an $\UF$-cover.

A multicovered space $(X,\mu)$ is said to be {\em $\mathsf F$-Menger} if for each sequence $\la U_n\ra_{n\in\w}\subset\mu$ of covers of $X$ there is an $\mathsf F$-cover $\la B_n\ra_{n\in\w}$ of $X$ by $u_n$-bounded subsets $B_n\subset X$, $n\in\w$. More information on $\mathsf F$-Menger multicovered spaces can be found in \cite{BZ06}.

\begin{proposition}\label{Menger-Scheepers} A multicovered space $X$ is Menger (resp. Scheepers) if and only if $X$ is $\SF$-Menger (resp. $\UF$-Menger).
\end{proposition}

\begin{proof} The ``if'' part trivially follows from the definitions.

To prove the ``only if'' part, assume that a multicovered space $(X,\mu)$ is Menger. To prove that it is $\SF$-Menger, fix a sequence of covers $\la u_n\ra_{n\in\w}\in\mu^\w$. Since $X$ is Menger, for every $k\in\w$ there is a cover $\la B_{n,k}\ra_{n\ge k}$ of $X$ by $u_n$-bounded subsets $B_{n,k}\subset X$. It follows that for every $n\in\w$ the finite union $B_n=\bigcup_{k\le n}B_{n,k}$ is $u_n$-bounded in $X$. We claim that $\la B_n\ra_{n\in\w}$ is an $\SF$-cover of $X$. Observe that for every point $x\in X$ and every $k\in\w$ the set $\{n\ge k:x\in B_{n,k}\}$ is not empty. Consequently, the set $\{n\in \w:x\in B_n\}\supset\bigcup_{k\in\w}\{n\ge k:x\in B_{n,k}\}$ is infinite and hence the family $\big\{\{n\in\w:x\in B_n\}:x\in X\big\}$ can be enlarged to a semifilter $\F\in\SF$, witnessing that $\la B_n\ra_{n\in\w}$ is an $\SF$-cover of $X$.

If the multicovered sapce $X$ is Scheepers, then the covers $\la B_{n,k}\ra_{n\in\w}$, $k\in\w$, can be chosen to be $\w$-covers. This means that for any finite subset $F\subset X$ and every $k\in\w$ the set $\{n\ge k:F\subset B_{n,k}\}$ is not empty and consequently, the set $\{n\in\w:F\subset B_n\}\supset\bigcup_{k\in\w}\{n\ge k:F\subset B_{n,k}\}$ is infinite. It follows that the family $\big\{\{n\in\w:x\in B_n\}:x\in X\big\}$ can be enlarged to a free ultrafilter $\F\in\UF$ witnessing that $\la B_n\ra_{n\in\w}$ is an $\UF$-cover of $X$.
\end{proof}

The $\mathsf F$-Menger property behaves nicely only for relatively large families $\mathsf F\subset\SF$, containing together with each semifilter $\F$ all its finite-to-one images. We recall that a function $\varphi:\w\to\w$ is called {\em finite-to-one} if for each $y\in\w$ the preimage $\varphi^{-1}(y)$ is finite and nonempty. It is clear that each monotone surjection $\varphi:\w\to\w$ is finite-to-one.

Given a family of semifilters $\mathsf F\subset\SF$ consider two its extensions:
$$\mathsf F_\asymp=\bigcup_{\F\in\mathsf F}[\F]\mbox{ \ and \ }\mathsf F_\da=\bigcup_{\F\in\mathsf F}\F_\da$$where $\F_\da=\{\varphi(\F)\mid\mbox{$\varphi:\w\to\w$ is a monotone surjection}\}$. It is clear that $$\mathsf F\subset\mathsf F_\da\subset\mathsf F_{\asymp}.$$

The following important lemma allows us to reduce the $\mathsf F$-Menger property of centered multicovered spaces to the $\F_\da$-Menger property for a suitable semifilter $\F\in\mathsf F$.

\begin{lemma} \label{count_cof} If a multicovered space $(X,\mu)$ with countable cofinality is $\mathsf F_\asymp$-Menger for some family of semifilters $\mathsf F\subset\SF$, then $(X,\mu)$ is $\F_\da$-Menger for some semifilter $\F\in\mathsf F$.
\end{lemma}

\begin{proof}
Let $\la u_n\ra_{n\in\w}\in\mu^\w$ be a witness for $\cf(X,\mu)=\w$ such that each $u_{n+1}$-bounded
subset of $X$ is $u_n$-bounded.
Since $(X,\mu)$ is $\mathsf F_\asymp$-Menger, there exists $\F\in\mathsf F$, a finite-to-finite
multifunction $\Phi:\w\Rightarrow\w$,  and a cover
$\la A_n\ra_{n\in\w}$ of $X$  such that each $A_n$ is $u_n$-bounded and $\Phi(\N_x)\in\F$,
where $\N_x=\{n\in\w: x\in A_n\}$, $x\in X$.

 Let $\la n_m\ra_{m\in\w}$ be an increasing
sequence of natural numbers such that $n_0=0$ and $n_{m+1}>1+\max(\Phi\cup\Phi^{-1})(\{0,\dots,n_m\})$ for $m\in\w$. It follows that
$n_{m+1}>m+1$ and $\Phi([n_{m+1},n_{m+2}))\subset [n_m,n_{m+3})$
for every $m\in\w$.

Set $C_m=\bigcup_{n\in [n_{m-1}, n_{m+3})}A_n$ and note that
it is  $u_m$-bounded.
Let $\phi:\w\to\w$ be the monotone surjection such that $\phi^{-1}(m)=[n_m, n_{m+1})$.
We claim that $  \{m\in\w:x\in C_m\}\in\phi(\F) $
for every $x\in X$. For this sake we shall show that $\{m\in\w:x\in C_m\}^\ast\supset\phi(\Phi(\N_x))$. Indeed, if $n\in\N_x\cap [n_{m+1}, n_{m+2})$
for some $m\in\w$, then $\phi(\Phi(n))\subset\{m,m+1, m+2\}$.
By the definition of $C_j$'s,
$A_n\subset\bigcup_{n\in [n_{m+1}, n_{m+2})}A_n\subset C_j$
for   every $j\in\{m,m+1, m+2\}$, which means that $\{m,m+1, m+2\}\subset \{m\in\w:x\in C_m\}$, and hence $\phi(\Phi(\N_x))\subset^\ast \{m\in\w:x\in C_m\}$.
Thus $\la C_m\ra_{m\in\w}$ is an $\{\phi(\F)\}$-cover of $X$ by $u_m$-bounded subsets.

For every $x\in X$ we denote by $\N'(x)$ the set $\{m\in\w:x\in C_m\}$.
Given $\la v_k\ra_{k\in\w}\in\mu^\w$, we can find an
increasing  sequence $\la m_k\ra_{k\in\w}$ of natural numbers such that each $u_{m_k}$-bounded
subset of $X$ is $v_k$-bounded.
Let $B_k=\bigcup_{m\in [m_k,m_{k+1})} C_m$
and $\psi:\w\to\w$ be such that $\psi^{-1}(k)= [m_k,m_{k+1}).$
By the definition of $m_k$'s, $B_k$ is $v_k$-bounded.
 It follows  from the above that
$$(\psi\circ\phi)(\F)\ni \psi(\N'_x)\subset^\ast \{k\in\w:x\in B_k\}$$
 for every $x\in X$, and hence
$\psi\circ\phi$
is a witness for $\la B_k:k\in\w\ra$, being an $\F_\da$-cover of $X$ by $v_k$-bounded
subsets, which completes our proof.
\end{proof}

\begin{corollary}\label{pofig} Let $\mathsf F\subset\SF$ be a family of semifilters. A centered multicovered space $X$ is $\mathsf F_\asymp$-Menger if and only if $X$ is $\mathsf F_\da$-Menger.
\end{corollary}
\begin{proof}
Observe that a multicovered space
$(X,\mu)$ is  $\mathsf F_\asymp$- (resp. $\mathsf F_\da$-) Menger
if and only if such is also $(X,\nu)$ for every
countable centered $\nu\subset\mu$.
We can now  apply  Lemma~\ref{count_cof}.
\end{proof}

Therefore for any semifilter $\F$ the $[\F]$-Menger and $\F_\da$-Menger properties of centered multicovered spaces are equivalent. Of these two properties, the $\F_\da$-Menger property is better for applications while the $[\F]$-Menger property is more convenient for proofs.

For any free filter $\F$ on $\w$ the class of $\F_\da$-Menger spaces has the following
nice inheritance properties:

\begin{theorem}\label{op} Let $\F$ be a free filter on $\w$. Then:
\begin{enumerate}
\item Each subspace of an $\F_\da$-Menger multicovered space is $\F_\da$-Menger.
\item A  multicovered space $X$ is $\F_\da$-Menger if it can be written as the countable union $X=\bigcup_{n\in\w}A_n$ of $\F_\da$-Menger subspaces $A_n\subset X$, $n\in\w$.
\item A multicovered space $Y$ is $\F_\da$-Menger if it is the image of an $\F_\da$-Menger multicovered space $X$ under a uniformly bounded surjective map $f:X\to Y$.
\item The direct product $X\times Y$ of two  $\F_\da$-Menger
multicovered spaces is $\F_\da$-Menger.
\item The meet $\mu_0\wedge \mu_1$ of two  $\F_\da$-Menger multicovers $\mu_0$ and $\mu_1$ on a set $X$ is  $\F_\da$-Menger.
\end{enumerate}
\end{theorem}

\begin{proof}
(1) The first assertion is obvious.
\smallskip

(2) To prove the second property,
denote by $\mu$ the underlying multicover of $X$ and
fix any sequence of covers $\la u_m\ra_{m\in\omega}\in\mu^\omega$.

For every $n\in\w$ find a cover
$\la B_{n,m}\ra_{m\in\w}$ of the $\F_\da$-Menger subspace $X_n\subset X$ by $u_m$-bounded subsets $B_{n,m}\subset X$ such that
$$\big\{\{m\in\w:x\in B_{n,m}\}:x\in X_n\big\}\in\F_n$$
for some semifilter $\F_n\in\F_\da$.

For every $m\in\w$ consider the $u_m$-bounded subset $B_m=\bigcup_{n\le m}B_{n,m}$ of $X$. We claim that
$$\big\{\{m\in\w:x\in B_m\}:x\in X\big\}\in\F_\infty$$
where $\F_\infty=\bigcup_{n\in\w}\F_n$. Indeed, given any $x\in X$ we can find $n\in\w$ such that $x\in X_n$ and conclude that
$$\F_\infty\supset\F_n\ni \{m\ge n:x\in B_{n,m}\}\subset\{m\in\w:x\in B_m\}$$and thus $\{m\in\w:x\in B_m\}\in\F_\infty$.

We claim that $\F_\infty\in[\F]$. It follows that
 for every $n\in\w$ there is a finite-to-finite multifunction $\Phi_n:\w\Ra\w$
with $\Phi_n(\F_n)\subset\F$. It is easy to construct a finite-to-finite
 multifunction $\Phi:\w\Ra\w$ such that for every $n\in\w$ the inclusion
 $\Phi_n(k)\subset \Phi(k)$ holds for all but finitely many numbers $k$.
 Then for every $n\in\w$ we get
$$\Phi(\F_n)\subset \Phi_n(\F_n)\subset\F$$ and thus $\Phi(\F_\infty)\subset\F$, witnessing $\F_\infty\Subset\F$.
The inverse relation $\F_\infty\Supset \F$ follows from $\F\asymp\F_0\subset\F_\infty$.

Therefore the space $X$ is $[\F]$-Menger and hence $\F_\da$-Menger by Corollary~\ref{pofig}.
\smallskip

(3) Assume that $X$ is an $\F_\da$-Menger multicovered space and let
 $f:X\to Y$ be a uniformly bounded map onto a multicovered space $Y$.
 To show that $Y$ is $\F_\da$-Menger, fix a sequence of covers
$\la u_n\ra_{n\in\w}\in\mu_Y^\w$. By the uniform boundedness of $f$
there is a sequence of covers $\la v_n\ra_{n\in\w}\in\mu_X^\w$ such that
for every $n\in\w$ the image $f(B)$ of any $v_n$-bounded subset
$B\subset X$ is $u_n$-bounded. Using the $\F_\da$-Menger property of
$X$, find an $\F_\da$-cover $\la B_n\ra_{n\in\w}$ of $X$ by $v_n$-bounded
subsets $B_n\subset X$. Then each set $f(B_n)$ is $u_n$-bounded
in $Y$. We claim that $\la f(B_n)\ra_{n\in\w}$ is an $\F_\da$-cover of
$Y$. Since $\la B_n\ra_{n\in\w}$ is an $\F_\da$-cover of $X$, there is a
semifilter $\F'\in\F_\da$ such that $\big\{\{n\in\w:x\in B_n\}:x\in X\big\}\in\F'$.

Fix a point $y\in Y$ and find any $x\in f^{-1}(y)$ (which exists
by the surjectivity of $f$). Observe that $x\in B_n$ implies
$y\in f(B_n)$ and thus $$\F'\ni \{n\in\w:x\in
B_n\}\subset\{n\in\w:y\in f(B_n)\}$$ and finally $\{n\in\w:y\in f(B_n)\}\in\F'\in\F_\da$, witnessing that $\la f(B_n)\ra_{n\in\w}$ is an $\F_\da$-cover of $Y$.
\smallskip

(4) Let $(X,\mu_X)$ and $(Y,\mu_Y)$ be two $\F_\da$-Menger multicovered spaces. Since the $\F_\da$-Menger and $[\F]$-Menger properties are equivalent, it suffices to check that the direct product $(X\times Y,\mu_X\boxtimes\mu_Y)$ is $[\F]$-Menger. Fix a sequence of covers $\la w_k\ra_{k\in\omega} \in(\mu_X\boxtimes\mu_Y)^\w$.
 For every $ k\in\omega $ we can write $ w_{k} $ in the
form $ w_k=u_{k}\boxtimes v_{k}$, where
$ u_{k}\in\mu_X $ and $v_k\in\mu_Y$.
Without loss of generality,
each $u_{k+1}$- (resp. $v_{k+1}$-) bounded subset is $u_k$- (resp. $v_k$-) bounded.
Using the $[\F]$-Menger  property of $X$ and $Y$, find  $[\F]$-covers
$\la A_k\ra_{k\in\w}$ and $\la B_k\ra_{k\in\w}$ of $X$ and $Y$ respectively
such that each $A_k$ (resp. $B_k$) is $u_k$-bounded (resp. $v_k$-bounded).
Find a finite-to-finite multifunction $\Phi$ such that
$$\Phi\big(\big\{\{n\in\w:x\in A_n\}:x\in X\big\}\cup \big\{\{n\in\w:y\in B_n\}:
y\in Y\big\}\big)\subset\F.$$
 Let $\la n_k\ra_{k\in\w}$ be an increasing
sequence of natural numbers such that $n_0=0$ and $n_{k+1}>\max(\Phi\cup\Phi^{-1}(\{0,\dots,n_k\})$ for $k\in\w$. It follows that
$\Phi([n_{k+1},n_{k+2}))\subset [n_k,n_{k+3})$
for every $k\in\w$.
Set $A'_k=\bigcup_{n\in [n_k,n_{k+3})}A_n$,
$B'_k=\bigcup_{n\in [n_k,n_{k+3})}B_n$, and $\Psi(k)=  [n_k,n_{k+3})$.
We claim that the multifunction $\Psi:\w\Ra\w$ is a witness  for
$\la A'_k\times B'_k\ra_{k\in\w}$ being an $[\F]$-cover of $X\times Y$.
Indeed, let us fix $(x,y)\in X\times Y$ and set $F_x=\{n\in\w:x\in A_n\}$,
$F_y=\{n\in\w:y\in B_n\}$, and $F'_{x,y}=\{k\in\w:(x,y)\in A'_k\times B'_k\}$.
 It suffices to show that $\Phi(F_x)\cap\Phi(F_y)\subset \Psi(F_{x,y})$.
Given any $m\in \Phi(F_x)\cap\Phi(F_y)$, find $n_x\in F_x$ and $n_y\in F_y$
such that $m\in\Phi(n_x)\cap\Phi(n_y)$. Let $k$ be such that $m\in [n_{k+1},n_{k+2})$.
By our choice of the sequence $\la n_k\ra_{k\in\w}$, $\{n_x,n_y\}\subset [n_k, n_{k+3})$,
and hence $x\in A_{n_x}\subset A'_k$ and $y\in B_{n_y}\subset B'_k$, which means
that $k\in F_{x,y}$. By the definition of $\Psi$, $m\in [n_{k+1},n_{k+2})\subset\Psi(k)$,
which completes our proof.
\smallskip

(5) The fifth property is a direct consequence of the previous ones.
Indeed, consider the diagonal $\Delta_X=\{(x,x):x\in X\}$ of the product $X\times X$ endowed with the induced multicover $\mu=\mu_1\boxtimes\mu_2\uhr \Delta_X$.
It is easy to see that the map $$f:(\Delta_X,\mu)\to (X,\mu_0\wedge\mu_1),\;\; f:(x,x)\mapsto x$$ is uniformly bounded. By the fourth property, the direct product
$(X\times X,\mu_0\boxtimes\mu_1)$ is $\F_\da$-Menger, and hence so is its subspace $(\Delta_X,\mu)$. Finally, the third property implies that
$(X,\mu_0\wedge\mu_1)$ is $\F_\da$-Menger as well.
\end{proof}

\begin{corollary} \label{v14}
Let $\mathsf F$ be a family of free filters and $(X,\mu)$ be an $\mathsf F_\da$-Menger
multicovered space. Then all finite powers of $X$ are $\mathsf F_\da$-Menger.
In particular, the Scheepers property is preserved by finite powers.
\end{corollary}

\begin{proof}
It suffices to prove that $(X^2,\mu\boxtimes \mu)$ is $\mathsf F_\da$-Menger (by induction this will imply that the powers $X^{2^n}$, $n\in\w$, all are $\mathsf F_\da$-Menger).
Given a sequence $\la w_n\ra_{n\in\w}\in (\mu\boxtimes\mu)^\w$, write
each cover $w_n$ in the form $u_{n}\boxtimes v_n$,
and set $\nu=\mathtt{cen}(\{u_{n},v_n:n\in\w\}).$ Applying Lemma~\ref{count_cof}
we conclude that $(X,\nu)$ is $[\F]$-Menger for some $\F\in\mathsf F$,
and hence so is the square $(X^2,\nu\boxtimes\nu)$ by Theorem~\ref{op}(4).
It follows that there exists an $[\F]$-cover $\la C_n\ra_{n\in\w}$ of $X^2$
 such that each $C_n$
is a $w_n$-bounded subset of $X^2$, which finishes our proof.

The last assertion follows from the equivalence of the Scheepers and $\UF$-Menger properties, proved in Proposition~\ref{Menger-Scheepers}.
\end{proof}

For ultrafilters $\F$ coherent to no $Q$-point, the fourth assertion of Theorem~\ref{op} can be generalized to semi-direct products of $\F_\da$-Menger multicovered spaces. For this we shall need a characterization of such ultrafilters $\F$ in terms of the left subcoherence.

Following \cite[10.2.1]{BZd}, we define a semifilter $\F$ to be {\em left subcoherent} to a semifilter $\U$ (and denote this by $\F\uss\U$) if for every monotone unbounded
 $f:\w\to\w$ there is a finite-to-finite multifunction
 $\Phi:\w\Ra\w$ such that $\Phi(\F)\subset\U$ and
 $\Phi(n)\subset [0, f(n)]$
for all $n\in\w$. By Proposition 10.2.2 \cite{BZd}, for two semifilters $\F,\U$
 the relations $\F\uss\F\Subset\U$ imply $\F\uss\U$.
By Proposition 10.2.3 of \cite{BZd}, an ultrafilter is not coherent to a $Q$-point if and only if \ $\U\uss\U$.

\begin{theorem}\label{semiprod} If a free ultrafilter $\F$ on $\w$ is not coherent to a $Q$-point, then for any centered $\F_\da$-Menger multicovered spaces $X,Y$ their semi-direct product $X\bowtie Y$ is $\F_\da$-Menger.
\end{theorem}

\begin{proof} By definition, the multicover $\mu_X{\bowtie}\mu_Y$ of the semidirect product $X\bowtie Y$ is the meet $$\mu_X{\bowtie}\mu_Y=(\mu_X{\ltimes}\mu_Y)\wedge(\mu_X{\rtimes}\mu_Y)$$ of the multicovers of the left and right semi-direct products of the multicovered spaces      $(X,\mu_X)$ and $(Y,\mu_Y)$. According to Theorem~\ref{op}(5) and Corollary~\ref{pofig}, the $\F_\da$-Menger property of the multicovered space $(X\times Y,\mu_X\bowtie\mu_Y)$ will follow as soon as we prove the $[\F]$-Menger property for the left and right semi-direct products $X\ltimes Y$ and $X\rtimes Y$.

To prove the $[\F]$-Menger property of $X\ltimes Y$, fix a sequence $\la w_n\ra_{n\in\w}\in(\mu_X\ltimes \mu_Y)^\w$
 of covers of $X\times Y$.
It follows from the definition of the multicover $\mu_X\ltimes\mu_Y$
that for every cover $w_n$ there is a cover $u_n\in\mu_X$ such that for every $u_n$-bounded subset
$B\subset X$ there is a cover $v\in\mu_Y$ such that for every $v$-bounded subset $D\subset Y$ the product
 $B\times D$ is $w_n$-bounded.

The $[\F]$-Menger property of the multicovered space $(X,\mu_X)$ yields an $[\F]$-cover $\la B_n\ra_{n\in\w}$ of $X$ by $u_n$-bounded subsets $B_n\subset Y$. For every $n\in\w$ find a cover $v_n\in\mu_Y$ such that for each $v_n$-bounded subset $D\subset Y$ the product $B_n\times D$ is $w_n$-bounded. Since $(Y,\mu_Y)$ is centered we can additionally assume that each $v_{n+1}$-bounded subset of $Y$ is $v_n$-bounded for $n\in\w$.

 The $[\F]$-Menger property
 of the muticovered space $(Y,\mu_Y)$ yields an $[\F]$-cover $\la C_n\ra_{n\in\w}$
 of  $Y$ by $v_n$-bounded subsets $C_n\subset Y$.
Since $\la B_n\ra_{n\in\w}$ and $\la C_n\ra_{n\in\w}$ are $[\F]$-covers, there are semifilters $\F_X,\F_Y\in[\F]$ such that
$$\big\{\{n\in\w:x\in B_n\}:x\in X\big\}\subset\F_X\mbox{ and }\big\{\{n\in\w:y\in C_n\}:y\in Y\big\}\subset\F_Y.$$

From $\F_X\Subset\F$ we conclude that there is a finite-to-finite multifunction $\Phi:\w\Ra\w$ such that $\Phi(\F_X)\subset\F$.
Since the coherence class $[\F]=[\F_Y]$ contains no $Q$-point, we can apply
Propositions 10.2.2 and 10.2.3 of \cite{BZd} to conclude that $\F_Y\uss\F$. Consequently, there is a finite-to-finite multifunction $\Psi:\w\Ra\w$ such that
$$\Psi(\F_Y)\subset\F\mbox{  and $\max\Psi(n)<\min \Phi([n,+\infty))$ for all $n\in\w$}.$$
It follows that $\Psi^{-1}\circ\Phi(n)\subset [n,\infty)$ for all $n\in\w$.

Let $D_n=\bigcup_{m\in\Psi^{-1}\circ \Phi(n)}C_m$ for every $n\in\w$. It follows from $\Psi^{-1}\circ \Phi(n)\subset [n,\infty)$ that each set $D_n$ is $v_n$-bounded. Consequently, the product $B_n\times D_n$ is $w_n$-bounded. We claim that $\la B_n\times D_n\ra_{n\in\w}$ is an $[\F]$-cover of $X\times Y$. Given any pair $(x,y)\in X\times Y$ it suffices to check that $\Phi(F_{(x,y)})\in\F$ where $F_{(x,y)}=\{n\in\w:(x,y)\in B_n\times D_n\}$.
Consider the sets $F_x=\{n\in\w:x\in B_n\}\in\F_X$ and $F_y=\{n\in\w:n\in C_n\}\in\F_Y$.

The inclusion $\Phi(F_{(x,y)})\in\F$ will follow as soon as we check that
$\Phi(F_x)\cap\Psi(F_y)\subset\Phi(F_{(x,y)})$. Take any number $m\in\Phi(F_x)\cap\Psi(F_y)$ and find points $n\in\Phi^{-1}(m)\cap F_x$ and $k\in\Psi^{-1}(m)\cap F_y$. It follows that $k\in\Psi^{-1}\circ\Phi(n)$ and thus $y\in C_k\subset D_n$. Consequently, $(x,y)\in B_n\times D_n$ and hence $n\in F_{(x,y)}$. Now we see that $m\in\Phi(n)\subset\Phi(F_{(x,y)})$  and thus $\Phi(F_x)\cap\Psi(F_y)\subset\Phi(F_{(x,y)})$. This completes the proof of the $[\F]$-Menger property of the left semi-direct product $X\ltimes Y$. By Corollary~\ref{pofig}, the space $X\ltimes Y$ has the $\F_\da$-Menger property.

The same property of the right semi-direct product $X\rtimes Y$ can be proved by analogy.
\end{proof}


\section{(Semi-)multicovered groupoids} \label{sem_mult_gr}

In this section we shall introduce a new mathematical object -- a semi-multicovered groupoid -- in which algebraic and multicover structures are connected.

By a {\em groupoid} we shall understand a set $X$ endowed with one binary operation $\cdot:X\times X\to X$ and one unary operation $(\cdot)^{-1}:X\to X$.

Note that a groupoid $X$ is a group if the binary operation is associative, $X$ has a two-sided unit $1$ and for every $x\in X$ we get $xx^{-1}=1=x^{-1}x$. However,
 the notion of a groupoid is much more general than one can expect. For example, each set
$X$ endowed with a binary operation $\cdot:X\times X\to X$ can be thought of
as a groupoid whose unary operation $(\cdot)^{-1}:X\to X$ is the identity.

A subset $A$ of a groupoid  $X$ is called a {\em sub-groupoid} of $X$ if $x\cdot y\in A$ and $x^{-1}\in A$ for all $x,y\in A$. We say that a groupoid $X$ is {\em algebraically generated\/} by a subset $A\subset X$ if $X$ coincides with the smallest sub-groupoid  of $X$ that contains the subset $A$.


By a {\em multicovered groupoid} (resp. {\em semi-multicovered groupoid\/}) we understand a groupoid $X$ endowed with a centered multicover $\mu$ such that the unary  operation $(\cdot)^{-1}:X\to X$ of $X$ is uniformly bounded and
the binary operation $\cdot:X\times X\to X$ of $X$ is uniformly bounded as a function from $(X\times X,\mu\boxtimes\mu)$ (resp. $(X\times X,\mu\bowtie\mu)$) to $(X,\mu)$.
It is clear that each multicovered groupoid is a semi-multicovered groupoid.

Topological groups endowed with the two-sided multicover are typical examples of semi-multicovered groupoids.

\begin{proposition} \label{good_ex}
A topological group $G$ endowed with the two-sided multicover $\mu_{L\wedge R}$ is a semi-multicovered groupoid. If $G$ is commutative, then $(G,\mu_{L\wedge R})$ is a multicovered groupoid.
\end{proposition}
\begin{proof}
The second (``commutative'') part is fairly easy. Therefore we give the proof
 only
of the first part.

The uniform boundedness of the inversion operation $(\cdot)^{-1}:G\to G$ trivially follows from the equality
$$(xU\cap Ux)^{-1}=U^{-1}x^{-1}\cap x^{-1}U^{-1}$$ holding for any point $x\in X$ and any neighborhood $U\subset G$ of the unit $1\in G$.

To prove the uniform boundedness of the binary operation $\cdot:G\bowtie G\to G$, we first prove that this operation is uniformly bounded as a map from $(G\times G,\mu_L\rtimes\mu_L)$ to $(G,\mu_L)$.

 Given any cover $u=\{xU:x\in G\}\in\mu_L$, find an open neighborhood $U_1$ of $1$
such that $U_1\cdot U_1\subset U$.  Given any $x\in G$,
find an open neighborhood
$V_x $ of $1$ such that $V_x \cdot x\subset U_1$. Set $u_1=\{x U_1: x\in G\}\in\mu_L$
and for every $x\in G$ set $v_x=\{y\cdot V_x :y\in G\}\in\mu_L$. It follows
from the above  that
if $A\subset G\times G$ is
$\{ yV_x\times xU_1: x,y\in G\}$-bounded,
then $\cdot (A)=\{a\cdot b: (a,b)\in A\}$ is $u$-bounded.
Since $u$ was chosen arbitrary,
 $\cdot : (G,\mu_L)\rtimes (G,\mu_L)\to (G,\mu_L)$ is uniformly bounded.
In the same way we can prove that
$\cdot: (G,\mu_R)\ltimes (G,\mu_R)\to (G,\mu_R)$ is uniformly bounded,
and hence $\cdot: G\times G\to G$ is uniformly bounded with respect
to the multicovers
$(\mu_L\rtimes\mu_L)\wedge (\mu_R\ltimes\mu_R)$ and $\mu_L\wedge\mu_R$.
Since the identity maps
$$ \id: (G\times G,\mu_{L\wedge R}\bowtie \mu_{L\wedge R})\to (G\times G, (\mu_R\rtimes\mu_R)\wedge (\mu_L\ltimes\mu_L)) $$
and $\id :\mu_L\wedge\mu_R \to \mu_{L\wedge R}$
are uniformly bounded
and the composition of two uniformly bounded maps is
uniformly bounded, the map
$$ \cdot : (G\times G,\mu_{L\wedge R}\bowtie \mu_{L\wedge R})\to
(G,\mu_{L\wedge R}) $$
is uniformly bounded as well, which finishes our proof.
\end{proof}


\section{The $\mathsf F$-Menger property in (semi-)multicovered groupoids}

In this section we characterize the $\mathsf F$-Menger property in (semi-)multicovered groupoids. We start with the $\w$-boundedness of  groupoids. For topological groups the following proposition was proved by I.~Guran \cite{Gu81}.

\begin{proposition}\label{ob} A semi-multicovered groupoid $X$ is $\w$-bounded if and only if it is algebraically generated by an $\w$-bounded subspace $A\subset X$.
\end{proposition}

\begin{proof} The ``only if'' part is trivial. To prove the ``if'' part, assume that $X$ is algebraically generated by an $\w$-bounded subspace $A\subset X$. Let $A_0=A$ and $A_{n+1}=A_n\cup A_n^{-1}\cup (A_n\cdot A_n)$ for $n\in\w$. By induction we shall show that each subspace $A_n\subset X$ is $\w$-bounded.

This is clear for $n=0$. Suppose that for some $n\in\w$ the subspace $A_n$ is $\w$-bounded.  The set $A_n^{-1}=\{a^{-1}:a\in A_n\}$ is $\w$-bounded, being the image of the $\w$-bounded multicovered space $A$ under the uniformly bounded map $(\cdot)^{-1}:X\to X$.
By Proposition~\ref{pob}, the semidirect product $A_n\bowtie A_n$ is  $\w$-bounded and so is its image $A_n\cdot A_n\subset X$ under the uniformly bounded binary operation $\cdot:X\bowtie X\to X$. Consequently, the subspace $A_{n+1}=A_n\cup A_n^{-1}\cup (A_n\cdot A_n)$ is $\w$-bounded, being the finite union of $\w$-bounded subspaces.

Now we see that the space $X$ is $\w$-bounded, being a
union of $\w$-bounded subspaces $A_n$, $n\in\w$.
\end{proof}

The following two theorems characterize the $\mathsf F$-Menger property in
(semi-) multi\-covered groupoids.

\begin{theorem} \label{bi1}
 Let $\mathsf F=\mathsf F_\da$ be a family of free filters on $\w$.
A multicovered groupoid $X$ is $\mathsf F$-Menger if and only if it
 is algebraically generated by an $\mathsf F$-Menger subspace $A\subset X$.
\end{theorem}

\begin{proof}
The ``only if'' part is obvious. To prove the ``if'' part, suppose that some  $\mathsf F$-Menger subspace $A$ of $X$ algebraically generates $X$. First assume that the multicovered space $X$ has countable cofinality. In this case we can apply Lemma~\ref{count_cof} to conclude that the multicovered space $A$ is $\F_\da$-Menger for some filter $\F\in\mathsf F$. Let $A_0=A$ and $A_{n+1}=A_n\cup A_n^{-1}\cup (A_n\cdot A_n)$ for $n\in\w$. By induction we shall show that for every $n\in\w$ the multicovered subspace $A_n$ of $X$ is $\F_\da$-Menger.

This is clear for $n=0$. Assuming that for some $n$ the space $A_n$ is $\F_\da$-Menger, we can apply Theorem~\ref{op}(3) and the uniform boundedness of the inversion $(\cdot)^{-1}:X\to X$ to conclude that the subset $A_n^{-1}$ of $X$ is $\F_\da$-Menger. Since the multicovered space $A_n$ is $\F_\da$-Menger, we can apply Theorem~\ref{op}(4) and conclude that the direct product $A_n\times A_n$ is $\F_\da$-Menger. Now the uniform boundedness of the multiplication $\cdot :X\times X\to X$ implies that the image $A_n\cdot A_n$ of $A_n\times A_n$ is $\F_\da$-Menger and so is the union $A_{n+1}=A_n\cup A_n^{-1}\cup (A_n\cdot A_n)$ according to Theorem~\ref{op}(2).

By Theorem~\ref{op}(2), the union $\bigcup_{n\in\w}A_n$ is $\F_\da$-Menger. The latter union coincides with $X$ because $X$ is algebraically generated by $A$. Therefore the multicovered space $X$ is $\F_\da$-Menger. Since $\F_\da\subset\mathsf F_\da=\mathsf F$, the space $X$ is $\mathsf F$-Menger.

Now we consider the general case (of arbitrary cofinality of $X$). Let $\mu$ be the multicover of $X$. To prove that $(X,\mu)$ is $\mathsf F$-Menger, take any sequence of covers $\langle u_n\rangle_{n\in\w}\in\mu^\w$. Using the uniform boundedness of the multiplication and inversion it is easy to find a countable centered subfamily $\mu'\subset\mu$ such that $\{u_n\}_{n\in\w}\subset\mu'$ and the operations $\cdot:(X\times X,\mu'\times\mu')\to (X,\mu')$, $(\cdot)^{-1}:(X,\mu')\to (X,\mu')$  are uniformly bounded.

Since the identity map $\id:(X,\mu)\to (X,\mu')$ is uniformly bounded, the
$\mathsf F$-Menger property of $(A,\mu{\uhr} A)$ implies the $\mathsf F$-property
of $(A,\mu'{\uhr} A)$. Since the multicovered space $(X,\mu')$ has countable cofinality, the preceding case guarantees that $(X,\mu')$ is $\mathsf F$-Menger because the multicovered groupoid $(X,\mu')$ is algebraically generated by the $\mathsf F$-Menger subspace $A\subset (X,\mu')$. Consequently, for the sequence $\la u_n\ra_{n\in\w}\in(\mu')^\w$ there is an $\mathsf F$-cover $\la B_n\ra_{n\in\w}$ by $u_n$-bounded subsets $B_n\subset X$. This
verifies that the multicovered space $(X,\mu)$ is $\mathsf F$-Menger.
\end{proof}

A similar result also holds  for semi-multicovered groupoids.

\begin{theorem} \label{bi2} Assume that a family $\mathsf F=\mathsf F_\da$ of free ultrafilters on $\w$ contains no $Q$-point.  A semi-multicovered groupoid $X$ is $\mathsf F$-Menger if $X$ is algebraically generated by an $\mathsf F$-Menger subspace $A\subset X$ and one of the following conditions holds:
\begin{enumerate}
\item $\mathsf F\subset[\F]$ for some $\F\in\mathsf F$;
\item The multicovered space $X$ has countable cofinality;
\item The multicovered space $X$ is properly $\w$-bounded;
\item The multicovered space $X$ is   uniformizable.
\end{enumerate}
\end{theorem}
\begin{proof}  (1) Assume that $\mathsf F\subset[\F]$ for some $\F\in\mathsf F$. Replacing the direct products in the proof of Theorem~\ref{bi1} by  semi direct products and applying Theorem~\ref{semiprod} instead of Theorem~\ref{op}(4) we can show that $X$ is $\mathsf F$-Menger.
\smallskip

(2) If the multicovered space $X$ has countable cofinality, then by Lemma~\ref{count_cof} the $\mathsf F$-Menger property of $A$ implies the $\F_\da$-Menger property of $A$ for some ultrafilter $\F\in\mathsf F$. By the first case, the multicovered space $X$ is $\F_\da$-Menger and consequently  $\mathsf F$-Menger.
\smallskip

(3) If the multicovered space $X$ is properly $\w$-bounded, then we can repeat the argument of the proof of Theorem~\ref{bi1}, replacing direct products
 with  semi-direct products everywhere and applying Theorem~\ref{semiprod} instead of
 Theorem~\ref{op}(4). The existence of $\mu'$,
which was straightforward for direct products, now follows from
Lemma~\ref{to_count}.
\smallskip

(4) Assume that the multicovered space $X$ is uniformizable. The subspace $A\subset X$, being $\mathsf F$-Menger, is $\w$-bounded. By Propositions~\ref{ob} and \ref{unif}, the uniformizable space $X$ is (properly) $\w$-bounded. Now the previous case completes the proof.
\end{proof}

Theorems~\ref{bi1}, \ref{bi2} and Proposition~\ref{Menger-Scheepers} imply the following characterization of the Scheepers property in \mbox{(semi-)multicovered} groupoids.

\begin{corollary} \label{cor_bi1} A multicovered groupoid $X$ is Scheepers if and only if $X$ is algebraically generated by a Scheepers subspace $A\subset X$.
\end{corollary}

\begin{corollary} \label{cor_bi2} If no $Q$-point exists, then a
semi-multicovered groupoid $X$ is Scheepers provided $X$ is algebraically
generated by a Scheepers subspace $A\subset X$ and one of the following conditions holds:
\begin{enumerate}
\item Any two ultrafilters are coherent;
\item The multicovered space $X$ has countable cofinality;
\item The multicovered space $X$ is properly $\w$-bounded;
\item The multicovered space $X$ is uniformizable.
\end{enumerate}
\end{corollary}

Note that the first condition
of Corollary~\ref{cor_bi2} is nothing else but the  NCF principle. Recall that
 no $Q$-points exist under NCF.


\section{The $\mathsf F$-Menger property in topological groups}\label{F-Menger-groups}

In this section we apply the general results proved in the preceeding section to studying the $\mathsf F$-Menger property in topological groups. According to Proposition~\ref{good_ex}, each topological group $G$ endowed with the two-sided multicover $\mu_{L\wedge R}$ is a semi-multicovered groupoid. Now we can
 apply Corollaries~\ref{cor_bi1} and \ref{cor_bi2} in order to prove:

\begin{corollary}\label{cor6_1}
 Let $\mathsf F=\mathsf F_\da$ be a family of free filters on $\w$. An  abelian
 topological group $G$ endowed with the two-sided multicover $\mu_{L\wedge R}$ is
 $\mathsf F$-Menger if and only if the group $G$ is algebraically generated by an
$\mathsf  F$-Menger subspace $A\subset G$.
\end{corollary}

\begin{corollary} \label{cor6_2} Assume that a family $\mathsf F=\mathsf F_\da$  of free ultrafilters on $\w$ contains no Q-points. A topological group $G$ endowed with the multicover $\mu_{L\wedge R}$ is $\mathsf F$-Menger if and only if the group $G$ is algebraically  generated by an $\mathsf F$-Menger subspace $A\subset G$.
\end{corollary}

\begin{corollary} \label{cor6_3} If no Q-point exists, then a topological group $G$ endowed with the multicover $\mu_{L\wedge R}$ is Scheepers if and only if the group $G$ is algebraically  generated by a Scheepers subspace $A\subset G$.
\end{corollary}

The $Q$-point assumption cannot be removed from Theorems~\ref{semiprod}, \ref{bi2} and Corollaries~\ref{cor_bi2}, \ref{cor6_2}, and \ref{cor6_3}.
To construct a suitable counterexample, consider the homeomorphism group $H(\IR_+)$ of the half-line $\IR_+=[0,\infty)$, endowed with the compact-open topology. This group will be considered as a semi-multicovered groupoid  endowed with the two-sided multicover $\mu_{L\wedge R}$.

\begin{theorem}\label{ess_q}
\begin{enumerate}
 \item[$(1)$] The semi-multicovered groupoid
$H(\IR_+)$ is not Menger and hence not Scheepers.
\item[$(2)$] For every Q-point $\F$ there is an $\LL_\da$-Menger subspace $A\subset(H(\IR_+),\mu_{L\wedge R})$ generating the group $H(\IR_+)$ in the sense that $H(\IR_+)=A\circ A\circ A\circ A $.
\end{enumerate}
\end{theorem}

\begin{proof}
We recall that the  compact-open topology on $H=H(\IR_+)$ is generated by the base
consisting of sets of the form
$$ [f,n]=\{g\in H:\forall x\leq n\: (|f(x)-g(x)|<2^{-n})\}, $$
where $f\in H$ and $n\in \w$.
Let $U_n=[\mathrm{id}_{\IR_+}, n]$. The family
 $\{U_n:n\in\w\}$ is a local base at the neutral element
$\mathrm{id}_{\IR_+}$ of $H$  and $U_{n+1}^2\subset
U_n$.
\smallskip

(1) The failure of the Menger property of the multicovered spaces $(H,\mu_{L\wedge R})$ and $(H,\mu_L)$ will follow as soon as we show that $H\ne\bigcup_{n\in\w}F_n\circ U_n$ for every sequence $\la F_n\ra_{n\in\w}$ of finite subsets of $H$. Given such a sequence, put $a_n=\max\{f(n+1):f\in F_n\}+1$.
For every $n\in\w$ and $(f,g)\in F_n\times U_n$ we have $f(g(n))\leq f(n+1)< a_n$.
Let $h\in H$ be the piecewise linear map such that $h(n)=a_n$.
It follows that $h\not\in\bigcup_{n\in\w}F_n\circ U_n$.
\smallskip

(2) In order to prove the second assertion it suffices to find
a sequence
$\la K_n\ra_{n\in\w}$ of finite subsets of $H$ such that the set
$$ (\bigcup_{L\in\LL}\bigcap_{n\in L} U_n K_n)\bigcap (\bigcup_{L\in\LL}\bigcap_{n\in L} U_n K_n)^{-1}$$
generates $H$.
 Let us consider a sequence $\la K_n:n\in\w\ra$
of finite subsets of $H$ such that $K_n= K_n^{-1}$ and
$W_n:=\{f\in H:\forall x\leq n\; (\frac{1}{3} x < f(x)<3x)\}\subset K_n U_n^{-1}\cap U_n K_n$
for all $n\in\w$. We claim that this sequence is as required. Indeed, let us note that
$$
(\bigcup_{L\in\LL}\bigcap_{n\in L} U_n K_n)\bigcap (\bigcup_{L\in\LL}\bigcap_{n\in L} U_n K_n)^{-1}=
\bigcup_{L\in\LL}\bigcap_{n\in L} U_n K_n\cap K_n U_n^{-1} \supset \bigcup_{L\in\LL}\bigcap_{n\in L} W_n.
$$
Let us fix $h\in H$ such that $h(x)>x$ for all $x>0$. Since $\LL$ is
 a $Q$-point, there exists
 $L=\{l_i:i\in\w\}\in\LL$ such  that $l_0>3$ and $l_{i+1}>2h(l_i+1)$ for all
 $i\in\w$.
  Now let $\phi$ be the
piecewise linear function such that $\phi(l_i)=3 l_i-1/2$ and
$\phi(l_i +1)=3l_{i+1}-1$ for all $i$. It is clear that $\phi\in
\bigcap_{n\in L} W_n$. Then $h(x)< (\phi\circ\phi)(x)$ for every $x>
0$. Indeed, fix $x> 0$. The following cases are possible.

a) $x\in [l_i, l_i +1)$ for some $i$. Then $\phi(x)\geq 3 l_i
-1/2>l_i+1$, hence $\phi(\phi(x))>\phi(l_i+1)\geq 3
l_{i+1}-1>l_{i+1}$, and consequently $\phi(\phi(x))>
h(l_{i}+1)>h(x)$.

b) $x\in [l_i+1, l_{i+1})$. Then $\phi(x)\geq 3
l_{i+1}-1/2>l_{i+1}+1$, and hence
$\phi(\phi(x))>\phi(l_{i+1}+1)=3l_{i+2}-1>l_{i+2}>h(l_{i+1})>h(x)$.

Thus $h(x)<\phi(\phi(x))$ for every $x$. Consequently
$\phi^{-1}(x)<(\phi^{-1}\circ h)(x)<\phi(x)$, which means that
$\phi^{-1}\circ h \in \bigcap_{n\in L} W_n$, and hence
$h=\phi\circ\phi^{-1}\circ h\in (\bigcap_{n\in L} W_n)\circ
(\bigcap_{n\in L} W_n)$. It suffices to note that the set of all
$h\in H$ such that $h(x)>x$ for all $x>0$ generates $H$.
\end{proof}


\section{Topological monoids and  (semi-)multicovered binoids}\label{top-monoids}

Looking at the results of the preceding section the reader can notice that among four canonical multicovers $\mu_L$, $\mu_R$, $\mu_{L\wedge R}$, $\mu_{L\vee R}$ on  topological groups, we distinguished one: the two-sided multicover $\mu_{L\wedge R}$. The reason is that the inversion operation is uniformly bounded with respect to this multicover. This is not true anymore for the multicovers $\mu_L$ and $\mu_R$. To treat topological groups endowed with those multicovers we will simply forget about the inversion operation and  think of groups as sets with a binary operation. Such algebraic structures are called {\em binoids}.

More precisely, a {\em binoid\/} is a set endowed with a binary operation $\cdot:X\times X\to X$. If this operation is associative, then $X$ is called a {\em semigroup}. A semigroup with a two-sided unit $1$ is called a {\em monoid}. It is clear that each group is a monoid. Each binoid $X$, endowed with the identity unary operation $$(\cdot)^{-1}:X\to X,\;\; (\cdot)^{-1}:x\mapsto x,$$ becomes a groupoid. So all the results about groupoids concern also binoids.

By a {\em multicovered binoid\/} (resp. {\em semi-multicovered binoid\/}) we understand a binoid $X$ endowed with a multicover $\mu$ making the binary operation $\cdot :X\times X\to (X,\mu)$ of $X$  uniformly bounded with respect to the multicover $\mu\boxtimes\mu$ (resp. $\mu\bowtie\mu$) on $X\times X$. Each (semi-)multicovered binoid can be thought of as a (semi-)multicovered groupoid with the  identity unary operation.

Let us now  return to topological groups and observe that they are examples of topological monoids. By a {\em topological monoid} we understand a monoid $X$ endowed with a topology $\tau$ making the binary operation $\cdot:X\times X\to X$ of $X$ continuous. The four multicovers considered earlier on topological groups can be equally defined on topological monoids.

Namely, for each topological monoid $(X,\tau)$ we can consider:
\begin{itemize}
\item the left multicover $\mu_L=\big\{\{xU:x\in X\}:1\in U\in\tau\big\}$;
\item the right multicover  $\mu_R=\big\{\{Ux:x\in X\}:1\in U\in\tau\big\}$;
\item the two-sided multicover  $\mu_{L\wedge R}=\big\{\{Ux\cap xU:x\in X\}:1\in U\in\tau\big\}$;
\item the R\"olke multicover  $\mu_{L\vee R}=\big\{\{UxU:x\in X\}:1\in U\in\tau\big\}$;
\item the multicover $\mu_L\wedge\mu_R$.
\end{itemize}
Let us observe that each left shift $l_a:X\to X$, $l_a:x\mapsto ax$, is uniformly bounded with respect to the left multicover $\mu_L$. Similarly, all the right shifts are uniformly bounded with respect to the right multicover $\mu_R$.

For each commutative topological monoid the multicovers $\mu_L$, $\mu_R$ and $\mu_{L\wedge R}$ coincide while $\mu_{L\vee R}$ is equivalent to $\mu_{L\wedge R}$ and to $\mu_L\wedge\mu_R$. For topological monoids the multicover $\mu_{L}\wedge\mu_R$ is more important than $\mu_{L\wedge R}$ (which is equivalent to $\mu_{L}\wedge \mu_R$ in topological groups).

If $(X,\tau)$ is a topological group then those multicovers are generated by suitable uniformities on $X$. The same is true for $\w$-bounded topological monoids.

\begin{proposition}\label{pr7_1}
 Let $X$ be a topological monoid endowed with a multicover
 $\mu\in\{\mu_L,\mu_R,\mu_L\wedge\mu_R,\mu_{L\wedge R},\mu_{L\vee R}\}$.
If the multicovered space $(X,\mu)$ is $\w$-bounded, then it is properly $\w$-bounded and hence is uniformizable.
\end{proposition}

\begin{proof} Assume that the multicovered space $(X,\mu_L)$ is $\w$-bounded. To prove that it is properly $\w$-bounded, fix any cover $u\in\mu_L$ and find a neighborhood $U\subset X$ of the unit $1\in X$ such that $u=\{xU:x\in X\}$.
By the continuity of the semigroup operation at $1$, there is a neighborhood $V\subset U$ of $1$ such that $VV\subset U$. Consider the cover $v=\{xV:x\in X\}\in\mu_L$. By the
$\w$-boundedness of the multicover $\mu_L$ there is a countable subset $C\subset X$ such that $X=CV$. It follows that $u'=\{cU:c\in C\}$ is a countable subcover of $u$. We claim that each $v$-bounded subset of $X$ is $u'$-bounded. It suffices to check that for every $x\in X$ the set $xV$ is $u'$-bounded. Since $CV=X$, there is a point $c\in C$ such that $x\in cV$. Then $xV\subset cVV\subset cU\in u'$. Thus the multicover $\mu_L$ is properly $\w$-bounded and, being centered, is uniformizable according to Proposition~\ref{unif}.

The proof of the fact that the  $\w$-boundedness of $\mu_R$ implies its
proper $\w$-bounded\-ness is completely analogous.

If the multicover $\mu_L\wedge\mu_R$ is $\w$-bounded, then the multicovers $\mu_L$ and $\mu_R$ are (properly) $\w$-bounded and so is their meet $\mu_L\wedge\mu_R$.

Concerning the multicovers $\mu_{L\wedge R}$ and $\mu_{L\vee R}$,
one just has to make a minor changes in the proof above. We demonstrate this on
$\mu_{L\wedge R}$.
Suppose that $\mu_{L\wedge R}$
is $\w$-bounded. To prove that it is properly $\w$-bounded, fix any cover $u\in\mu_{L\wedge R}$ and find a neighborhood $U\subset X$ of the unit $1\in X$ such that $u=\{xU\cap Ux:x\in X\}$. By the continuity of the semigroup operation at $1$, there is a neighborhood $V\subset U$ of $1$ such that $VV\subset U$. Consider the cover $v=\{xV\cap Vx:x\in X\}\in\mu_{L\wedge R}$. By the
$\w$-boundedness of the multicover $\mu_{L\wedge R}$ there is a countable subset $C\subset X$ such that $X=\bigcup_{c\in C}cV\cap Vc$. It follows
that $u'=\{cU\cap Uc:c\in C\}$ is a countable subcover of $u$. We claim that each $v$-bounded subset of $X$ is $u'$-bounded. It suffices to check that for every
$x\in X$ the set $xV\cap Vx$ is $u'$-bounded. Since $X=\bigcup_{c\in C}cV\cap Vc$,
 there is a point $c\in C$ such that $x\in cV\cap Vc$. Then
$xV\subset cVV\subset cU\in u'$ and $Vx\subset VVc\subset Uc\in u'$,
and consequently $xV\cap Vx\subset cU\cap Uc$, and hence is $u'$-bounded.
\end{proof}

Now we shall detect multicovers on topological monoids turning them into\break \mbox{(semi-)multicovered} binoids.

\begin{proposition} \label{com_mo}
Each commutative topological monoid $X$ endowed with any of the equivalent multicovers $\mu_L$, $\mu_R$, $\mu_L\wedge\mu_R$, $\mu_{L\wedge R}$ or $\mu_{L\vee R}$ is a multicovered binoid.
\end{proposition}

\begin{proof}
The commutativity of $X$ easily implies that all of these multicovers are equivalent.
Therefore it is enough to check that $(X,\mu_L)$ is a multicovered binoid.
Given $u\in\mu_L$, write it in the form $\{xU:x\in X\}$ for some open neighborhood
$U$ of $1$, and find $V\ni 1$ such that $V\cdot V\subset U$.
Set $v=\{xV:x\in X\}\in\mu_L$. In order to prove that the operation
$\cdot : (X\times X,\mu_L\boxtimes\mu_L)\to (X,\mu_L)$ is uniformly bounded
we need to show that $(xV)\cdot (yV)$ is $u$-bounded for arbitrary
$x,y\in X$. Since $X$ is commutative, $xVyV=xyVV\subset xyU\in u$, which finishes
our proof.
\end{proof}

For non-commutative topological monoids the situation is a bit more complicated.
Let us define a point $x\in X$ of a topological monoid $X$ to be {\em left balanced} (resp. {\em right balanced}) if for every neighborhood $U\subset X$ of the unit $1$ of $X$ there is a neighborhood $V\subset X$ of $1$ such that $Vx\subset xU$ (resp. $xV\subset Ux$). Observe that $x$ is left balanced if the left shift $l_x:X\to X$, $l_x:y\mapsto xy$, is open at $1$.

Let $B_{L}$ and $B_R$ denote respectively the sets of all left and right balanced  points of the monoid $X$.

A topological monoid $X$ is defined to be {\em left balanced} (resp. {\em right balanced}) if $X= B_L\cdot U$ (resp. $X= U\cdot B_R$) for every neighborhood $U\subset X$ of the unit $1$ in $X$. If a topological monoid $X$ is left and right balanced, then we say that $X$ is {\em balanced}.

Observe that the class of balanced topological monoids includes all commutative topological monoids.

\begin{proposition} \label{rb} If a topological monoid $X$ is balanced (resp. left balanced, right balanced), then $X$ endowed with the multicover $\mu_{L}\wedge \mu_R$ (resp. $\mu_L$, $\mu_R$)  is a semi-multicovered binoid.
\end{proposition}

\begin{proof} Assume that a topological monoid $X$ is left balanced and let $B_L$ be the set of left balanced points of $X$.
We shall show that the semigroup operation
$\cdot:(X\times X,\mu_L\rtimes \mu_L)\to (X,\mu_L)$ is uniformly bounded.

Given any cover $u\in\mu_L$, find a neighborhood $U_0\subset X$ of the unit $1$ of $X$ such that $u=\{xU_0:x\in X\}$. By the continuity of the operation at 1, there is a neighborhood $W\subset X$ of $1$ such that $WWW\subset U_0$. The monoid $X$, being left balanced,
 is equals to $B_L\cdot W$.
So, for every $y\in X$ we can find a left balanced point $b_y\in B_L$ such
that $y\in b_yW$. The left balanced property of
$b_y$ allows us to find a neighborhood $W_y\subset X$ of
$1$ such that $W_yb_y\subset b_yW$.

Now consider the cover $v=\{yW:y\in X\}$ and for every set $V\in v$ find a point
$y\in X$ with $V=yW$ and consider the cover $u_V=\{xW_y:x\in X\}\in\mu_L$.
Set $w=\{U\times V:V\in v,\; U\in u_V\}\in \mu_L\rtimes \mu_L$.
The uniform boundedness of the operation
$\cdot:(X\times X,\mu_L\rtimes\mu_L)\to(X,\mu_L)$ will follow as soon as we show that for every set $U\times V\in w$ the set $U\cdot V$ is $u$-bounded. Find $y\in X$ with $V=yW$ and $x\in X$ with $U=xW_y$. Now we see that $$U\cdot V=xW_yyW\subset xW_yb_yWW\subset xb_yWWW\subset xb_yU_0\in u.$$

By analogy we can prove that if the topological monoid $X$ is right balanced, then $X$ endowed with the multicover $\mu_R$ is a semi-multicovered binoid.
\smallskip

Now assume that the topological monoid $X$ is balanced. Since $X$ is both left and right balanced, the maps $\cdot:(X\times X,\mu_L\bowtie \mu_L)\to (X,\mu_L)$ and
$\cdot:(X\times X,\mu_R\bowtie \mu_R)\to (X,\mu_R)$ are uniformly bounded,
 and hence  so is the map
$$\cdot:\big(X\times X,(\mu_L\bowtie \mu_L)\wedge(\mu_R\bowtie \mu_R)\big)\to(X,\mu_L\wedge\mu_R).$$
Taking into account that the identity map
$$\id:\big(X\times X,(\mu_L\wedge \mu_R)\bowtie (\mu_L\wedge\mu_R)\big)\to \big(X\times X,(\mu_L\bowtie \mu_L)\wedge(\mu_R\bowtie \mu_R)\big)$$is uniformly bounded, we see that so is the map
$$\cdot:\big(X\times X,(\mu_L\wedge \mu_R)\bowtie (\mu_L\wedge\mu_R)\big)\to(X,\mu_L\wedge\mu_R).$$
\end{proof}

A typical example of a topological monoid is the semigroup $\w^\w$ of all maps
from  $\w$ to $\w$ endowed with the Tychonov product topology.
The semigroup operation on $\w^\w$ is given by the composition of functions.

Besides the five multicovers generated by the algebraic structure,
the monoid $\w^\w$ carries the product multicover $\mu_p$ on $\w^\w$. It consists
of the uniform covers
$u_n=\{[s]:s\in\w^n\}$, $n\in\w$, where $[s]=\{y\in\w^\w:y\uhr n =s\}$
(we identify $n$ with $\{0,\dots,n-1\}$ and write $f[n]$ for $\{f(0),\ldots,f(n-1)\}$).
Observe that the multicover $\mu_p$ coincides with the multicover $\mu_\rho$ generated by the
(standard) metric $\rho(x,y)=\inf\{2^{-n}:x(n)=y(n)\}$ on $\w^\w$.

The following proposition characterizes left and right balanced points in the topological monoid $\w^\w$.
We define a function $f:\w\to\w$ to be {\em eventually injective} if the restriction $f\uhr(\w\setminus n)$ is injective for some $n\in\w$.

\begin{proposition}\label{balanced} An element $f\in\w^\w$ of the topological monoid $\w^\w$ is
\begin{enumerate}
\item left balanced  if and only if the function $f$ is bounded or surjective;
\item right balanced  if and only if  $f$ is constant or eventually injective.
\end{enumerate}
\end{proposition}

\begin{proof} (1) To prove the ``if'' part of the first assertion, consider two cases.

1a) Suppose that $f$ is bounded and  fix $m\in\w$.
Set $M=\max\obr(f)+1$. We claim that $[\id\uhr M]\circ f\subset
f\circ [\id\uhr m]$. Indeed, a direct verification shows that
$[\id\uhr M]\circ f=\{f\}\subset f\circ [\id\uhr m]$.

1b) $f$ is surjective. Let us fix $m\in\w$. We claim that
$[\id\uhr f(m)]\circ f\subset f\circ [\id\uhr m]$. Indeed, fix $g\in [\id\uhr f(m)]$,
for every $k\geq m$ find $h(k)$ such that $f(h(k))=g(f(k))$ (it exists
by the surjectivity of $f$), and set $h\uhr m=\id\uhr m$. Then
$g\circ f=f\circ h$ and $h\in [\id\uhr m]$, and consequently
$g\circ f\in f\circ [\id\uhr m]$, which means that $f$ is left balanced.

To prove the ``only if'' part suppose that $f$ is unbounded and
there exists $p\in\w\setminus\obr(f)$.
We claim that $[\id\uhr M]\circ f\not\subset f\circ\w^\w$
for all $M\in\w$. Indeed, given $M$ find $n\in\w$ such that $f(n)\geq M$
and set $g(f(n))=p$ and $g\uhr(\w\setminus\{f(n)\})=\id$. Then
$p\in\obr(g\circ f)$ but $p\not\in\obr(f\circ h)$ for all $h$,
and consequently $g\circ f\not\in f\circ\w^\w$, which finishes our proof.
\medskip

(2) To prove the ``if'' part of the second assertion consider two cases:

2a) $f\in\w^\w$ is constant, i.e. there exists $n_0\in\w$ such that
$f(n)=n_0$ for all $n$. Let us fix $m\in\w$ and set $M=m$.
We claim that $f\circ [\id\uhr M]\subset [\id\uhr m]\circ f$.
Indeed, for every $g\in [\id\uhr M]$ and $n\in\w$ we have
$(f\circ g)(n)=n_0 = f(n)$, and hence $f\circ g=f = \id\circ f\in [\id\uhr m]\circ f$,
and consequently $f$ is left balanced.

2b) $f$ is eventually injective.
Let us fix $m\in\w$ and find $M\geq m$
such that
$$ (m\cup f(M))\cap f([M,+\infty))=\emptyset.   $$
We claim that $f\circ [\id\uhr M]\subset [\id\uhr m]\circ f$.
Fix $g\in [\id\uhr M]$ and for every $k\geq M$ set $h(f(k))=f(g(k))$.
Since $f\uhr [M,+\infty)$ is injective, this well defines a
 function $h: f([M,+\infty))\to\w$. Extend $h$ to the function on $\w$
by letting $h(l)=l$ for all $l\not\in f([M,+\infty))$.
The choice of $M$ guarantees that $h\in [\id\uhr m]$.
A direct verification shows that $f\circ g=h\circ f$,
which finishes the proof that $f$ is right balanced.

To prove the ``only if'' part suppose that $f$ is not constant and
for every $n\in\w$ there are distinct $m,l\geq n$ such that $f(l)=f(m)$.
Let us fix $p,q\in\w$ such that $f(p)\neq f(q)$.
We claim that $f\circ [\id\uhr M]\not\subset \w^\w\circ f$ for all $M\in\w$.
Given  $M\in\w$, find distinct $r,l\geq M$ such that $f(l)=f(r)$ and set
$g(r)=p, g(l)=q$, and $g(n)=n$ otherwise.
Then $(f\circ g)(r)\neq (f\circ g)(l)$, while $(h\circ f)(r)=(h\circ f)(l)$ for all $h\in\w^\w$,
and consequently $f\circ g\not\in\w^\w\circ f$, which finishes the proof.
\end{proof}

We are now  able to prove that the topological monoid $\w^\w$ is balanced.
 We recall that a multicovered
space $(X,\mu)$ is called {\em totally bounded} if $X$ is $u$-bounded for every cover
 $u\in\mu$. In this case we also say that the multicover $\mu$ is totally bounded. It is clear
 that any two totally bounded multicovers on a set $X$ are equivalent.

\begin{proposition} \label{od_mon}
\begin{enumerate}
\item The topological monoid $\w^\w$ is balanced.
\item The multicover $\mu_R$ is totally bounded.
\item The multicovers $\mu_L$ and $\mu_L\wedge\mu_R$ are equivalent
to the product multicover $\mu_p$.
\end{enumerate}
\end{proposition}

\begin{proof}
(1) Fix $m\in\w$ and  $x\in\w^\w$.

In order to prove that $\w^\w$ is left balanced
 we need to find a left balanced $y\in\w^\w$ and $g\in [\id\uhr m]$
such that $x= y\circ g$. Let $g$ be any injection
with $|\w\setminus\obr(g)|=\w$ and $g\uhr m=\id\uhr m$.
Define $y\uhr\obr g$ by $y(g(k))=x(k)$ (the correctness follows from the
injectivity of $g$) and extend $y$ onto $\w$ in such a way that
$\obr(y\uhr \w\setminus\obr(g))=\w $. Being surjective, $y$ is a left
 balanced according to Proposition~\ref{balanced},
 $g\in [\id\uhr m]$, and $x=y\circ g$.

Next, we prove that $\w^\w$ is right balanced.
We need to find a right balanced $y\in\w^\w$ and $g\in [\id\uhr m]$
such that $x= g\circ y$.
Set $y\uhr m=x\uhr m$, $F=m\cup x[m]$,
and   $g\uhr F =\id$. Let
 $y\uhr [m, +\infty)$ be an injection into $\w\setminus F$
and $g\uhr\w\setminus F$ be such that
 $g(y(k))=x(k)$ for all $k\geq m$.
It follows that
  $y$ is eventually
injective, $g\in [\id\uhr m]$, and $x=g\circ y$. By Proposition~\ref{balanced}, $y$ is right balanced.
\smallskip

(2) In fact, for every $n$ there exists $f\in\w^\w$
such that $[\id\uhr n]\circ f=\w^\w$. Indeed, any injection from $\w$
into $[n, +\infty)$ is obviously as required.
\smallskip

(3) Since $x\circ [\id\uhr m]\subset [x\uhr m]\in u_m$,
the identity map $\id:(\w^\w,\mu_L)\to(\w^\w,\mu_p)$ is uniformly bounded.
To show that
$\id: (\w^\w,\mu_p)\to (\w^\w,\mu_L)$ is uniformly bounded we need to prove
that for every $n\in\w$ there exists $m\in\w$ such that $[s]$ is
$\{x\circ [\id\uhr n]:x\in\w^\w\}$-bounded for all $s\in\w^m$.
We claim that $m=n$ is as required. Indeed, let us fix $s\in\w^n$
and define $x\in\w^\w$ letting $x\uhr n=s\uhr n$ and
extending it onto $\w$ in such a way that $\obr(y\uhr[n,+\infty))=\w$.
Given $y\in [s]$, for every $k\geq m$ find $h(k)$ such that $y(k)=x(h(k))$
(this is possible by the definition of $x$) and set $h\uhr m=\id\uhr m$.
It follows that $y=x\circ h\in x\circ [\id\uhr n]$. Since $y$ was chosen arbitrary,
we conclude that $[s]\subset x\circ [\id\uhr n]$, which finishes our proof
of the equivalence of $\mu_p$ and $\mu_L$.

Since $\mu_R$ is totally bounded, $\mu_L\wedge\mu_R$ is equivalent to $\mu_L$
and hence to $\mu_p$ as well.
\end{proof}

Other natural examples of balanced topological monoids arise from paratopological groups.
We recall that a {\em paratopological group} is a group $G$ endowed with a topology $\tau$ making the group operation continuous. The following statement is a direct consequence of Propositions~\ref{com_mo} and \ref{rb}.

\begin{corollary} \label{imp_ex}
\begin{enumerate}
\item Each abelian paratopological group $G$ endowed with the multicover $\mu_L\wedge \mu_R$ is a multicovered binoid.
\item Each paratopological group $G$ is a balanced topological monoid and endowed with one of the multicovers $\mu_L$, $\mu_R$ or $\mu_{L}\wedge \mu_R$ is a semi-multicovered binoid.
\end{enumerate}
\end{corollary}


\section{The $\mathsf F$-Menger property in topological monoids} \label{top_mon_s}

In this section we shall characterize the $\mathsf F$-Menger property in topological monoids. Combining Proposition~\ref{com_mo} with Theorem~\ref{bi1} we get the following corollaries:

\begin{corollary} \label{com_f} Let $\mathsf F=\mathsf F_\da$ be a family of free
filters on $\w$.
A commutative topological monoid $X$ endowed with the multicover $\mu_{L}\wedge \mu_R$ is
$\mathsf F$-Menger if and only if $X$ is algebraically
generated by an $\mathsf F$-Menger subspace $A\subset X$.
\end{corollary}

\begin{corollary} \label{co8_2}
A commutative topological monoid $X$ endowed with the multicover $\mu_{L}\wedge \mu_R$ is Scheepers if and only if $X$ is algebraically generated by a Scheepers subspace $A\subset X$.
\end{corollary}

The analogous results for non-commutative topological monoids  follow from
Propositions~\ref{ob}, \ref{pr7_1}, and \ref{rb} and Theorem~\ref{bi2}(3).

\begin{corollary} \label{co8_3} Assume that a family $\mathsf F=\mathsf F_\da$ of free ultrafilters on $\w$ contains no Q-point. A balanced (resp. left balanced, right balanced) topological monoid $X$ endowed with the multicover $\mu_{L}\wedge \mu_R$ (resp. $\mu_L$, $\mu_R$) is $\mathsf F$-Menger if and only if the monoid $X$ is algebraically generated by an $\mathsf F$-Menger subspace $A\subset X$.
\end{corollary}

\begin{corollary} \label{co8_4}
If no Q-point exists, then a balanced (resp. left balanced, right balanced) topological monoid $X$ endowed with the multicover $\mu_{L}\wedge \mu_R$ (resp. $\mu_L$, $\mu_R$) is Scheepers if and only if $X$ is algebraically generated by a Scheepers subspace $A\subset X$.
\end{corollary}

\section{Characterizing the $\mathsf F$-Menger property in free topological groups}
\label{letztlich}

We are now in a position  to present the proof of Theorem~\ref{main}, which is a special
 case of Corollary~\ref{cor9_2}.

We define a topological space $X$ to be $\mathsf F$-Menger if such is the multicovered space $(X,\mu_{\mathcal O})$ where $\mu_{\mathcal O}$ is the multicover consisting of all open covers of $X$. If the space $X$ is Lindel\"of, then the multicover $\mu_{\mathcal O}$ is equivalent to the multicover  $\mu_\U$ consisting of the uniform covers with respect to the universal uniformity of $X$, see
\cite[Corollary~15]{Zd06}.
The multicover $\mu_\U$ is equivalent to the multicover consisting of the covers $\{B_d(x):x\in X\}$ by 1-balls with respect to all continuous pseudometrics $d$ on $X$.

\begin{theorem} \label{main_techical}
Assume that a family $\mathsf F=\mathsf F_\da$ of free ultrafilters on $\w$ contains no Q-point. For a Tychonov space $X$ the  following conditions are equivalent:
\begin{enumerate}
\item All continuous metrizable images of $X$ are $\mathsf F$-Menger;
\item All continuous Lindel\"of regular images of $X$ are $\mathsf F$-Menger;
\item  The multicovered space $(X,\mu_\U)$ is $\mathsf F$-Menger;
\item  $(F(X),\mu_{L\wedge R})$ is $\mathsf F$-Menger;
\item  $(F(X),\mu_L)$ is $\mathsf F$-Menger;
\item   $(F(X),\mu_R)$ is $\mathsf F$-Menger;
\item   $(F(X),\mu_{L\vee R})$ is $\mathsf F$-Menger;
\item   $(A(X),\mu_{L\wedge R})$ is $\mathsf F$-Menger.
\end{enumerate}
\end{theorem}

\begin{proof}
The implications
$(4)\Ra (5,6)\Ra(7)\Ra(8)$
are straightforward.
\smallskip

The implication $(8)\Ra (3)$ follows from the fact that
 the restriction to $X$ of the natural uniformity of
$A(X)$ coincides with $\U_X$ \cite{Pe82}
and the  $\mathsf F_\downarrow$-Mengerness is   preserved by
taking subspaces, see Theorem~\ref{op}(1).
\smallskip

$(3)\Ra(2)$. Let $f:X\to T$ be a continuous map onto a Lindel\"of space
$T$. Then $f$ is uniformly continuous with respect to the universal uniformities $\U_X$ and $\U_T$,
and hence uniformly bounded with respect to the multicovers
$\mu_{\U_X}$ and $\mu_{\U_T}$, and consequently the multicovered space
$(T,\mu_{\U_T})$ is $\mathsf F$-Menger by Theorem~\ref{op}(3).
Applying \cite[Corollary~15]{Zd06}, we conclude that
the multicovers $\mu_{\mathcal O_T}$ and $\mu_{\U_T}$ of $T$ are equivalent,
and consequently $T$ is $\mathsf F$-Menger as a topological space.
\smallskip

The implication $(2)\Ra(1)$ is trivial.
\smallskip

$(1)\Ra (3)$. Suppose that $(X,\mu_\U)$ is not $\mathsf F$-Menger and fix a
sequence $\la U_n\ra_{n\in\w}\in\U^\w$ such that
$\la U_n(K_n)\ra_{n\in\w}$ is an $\mathsf F$-cover of $X$
 for no sequence $\la K_n\ra_{n\in\w}$ of finite subsets of $X$.
Let $d$ be a continuous pseudometric on $X$ such that
$\{(x,y): d(x,y)<1/n\}\subset U_n$ and $Y$ be the quotient-space of $X$
with respect to the equivalence relation $x\equiv y\leftrightarrow d(x,y)=0$.
A direct verification shows that $Y$ is a metrizable image of $X$ which fails to be
$\mathsf F$-Menger.
\smallskip

$(3)\Ra(4)$. Assume that the multicovered space $(X,\mu_\U)$ is $\mathsf F$-Menger. Consider the canonical embedding $i:X\to F(X)$ and observe that it is uniformly bounded as a map from $(X,\mu_\U)$ into $(F(X),\mu_{L\wedge R})$. By Theorem~\ref{op}(3), the image $i(X)\subset F(X)$ is $\mathsf F$-Menger and so is the semi-multicovered group $(F(X),\mu_{L\wedge R})$ according to Corollary~\ref{cor6_2}
\end{proof}

\begin{corollary} \label{cor9_2}
Assume that there is no $Q$-point. For a Tychonov space $X$ the  following conditions are equivalent:
\begin{enumerate}
\item All continuous metrizable images of $X$ are Scheepers;
\item All continuous Lindel\"of regular images of $X$ are Scheepers;
\item  The multicovered space $(X,\mu_\U)$ is Scheepers;
\item  $(F(X),\mu_{L\wedge R})$ is Scheepers;
\item   $(F(X),\mu_{L\vee R})$ is Scheepers;
\item   $(A(X),\mu_{L\wedge R})$ is Scheepers;
\item all finite powers of $F(X)$ are $o$-bounded;
\item   $F(X)$ is $o$-bounded;
\item   $A(X)$ is $o$-bounded.
\end{enumerate}
\end{corollary}
\begin{proof}
 The equivalence of the conditions $(1)$--$(6)$ is a special case
of Theorem~\ref{main_techical} for $\mathsf F=\mathsf{UF}$.
The implications $(7)\Rightarrow(8)$ and
$(8)\Rightarrow(9)$ are obvious.
Concerning $(9)\Rightarrow (1)$, it follows from Theorem~\ref{ab_ch}. Finally, the implication
$(4)\Rightarrow (7)$ follows from Proposition~\ref{v14}.
\end{proof}

\begin{remark} In spite of the example from Theorem~\ref{ess_q} we do not know if Theorem~\ref{main_techical} is true for all families $\mathsf F=\mathsf F_\da$ of ultrafilters (in particular, those containing $Q$-points).
\end{remark}

\end{document}